\documentclass[A4]{amsart}
\newtheorem{lemma}{Lemma}[subsection]
\newtheorem{propos}[lemma]{Proposition}

\newtheorem{theorem}[lemma]{Theorem}

\newtheorem{defin}[lemma]{Definition}

\newcommand{\CL}{\hbox{{$\mathcal L$}}}

\newcommand{\CT}{\hbox{{$\mathcal T$}}}

\newcommand{\cg}{\mathfrak{g}}

\newcommand{\C}{\mathbb{C}}
\newcommand{\R}{\mathbb{R}}

\newcommand{\diff}{{\rm diff}}

\newcommand{\g}{{\frak g}}

\newcommand{\del}{\partial}

\newcommand{\extd}{{\rm d}}

\newcommand{\isom}{{\cong}}

\newcommand{\tens}{\mathop{\otimes}}
\newcommand{\la}{{\triangleright}}
\newcommand{\ra}{{\triangleleft}}

\newcommand{\ad}{{\rm ad}}
\newcommand{\id}{{\rm id}}
\newcommand{\<}{\langle}
\renewcommand{\>}{\rangle}

\newcommand{\Tr}{{\rm Tr}}

\def\rcross{{\triangleright\!\!\!<}}
\def\lcross{{>\!\!\!\triangleleft}}

\newcommand{\varnabla}{{\hat \nabla}}

\newcommand{\eproof}{$\quad \heartsuit$\bigskip}
\newcommand{\nquad}{\kern-60pt}

\newcommand{\eqn}[2]{\begin{equation}#2\label{#1}\end{equation}}

\begin{document}

\title[Quantization by cochain twists and nonassociative differentials]{\rm\large
Quantization by cochain twists and nonassociative differentials}

\author{E.J. Beggs \& S. Majid}%

\address{EJB: Department of Mathematics\\University of Wales,
   Swansea SA2 8PP, UK\\
SM: School of Mathematical Sciences\\
Queen Mary, University of London\\ 327 Mile End Rd,  London E1
4NS, UK}



\date{18 June 2005/Rev. 15 Nov 2005}

\maketitle

\begin{abstract}We show that several  standard associative
quantizations in mathematical physics  can be expressed as cochain module-algebra twists  in the spirit
of Moyal products at least to $O(\hbar^3)$,  but to achieve this  we twist not by a 2-cocycle but by a 2-cochain. This implies a hidden nonassociavitity not visible in the algebra  itself but present in its deeper noncommutative differential geometry, a phenomenon first seen in our previous work on
semiclassicalisation of differential structures.   The quantisations are induced by  a classical group covariance and include: enveloping
algebras $U_\hbar(g)$ as quantisations of $g^*$, a Fedosov-type quantisation
of the sphere $S^2$ under a Lorentz group covariance, the Mackey quantisation of homogeneous spaces, and the standard quantum groups $\C_q[G]$. 
We also  consider the differential quantisation  of  $\R^n$ for a given symplectic
connection as part of our semiclassical analysis and we outline a proposal for the Dirac operator.
  \end{abstract}

\section{Introduction}

This paper is a sequel to \cite{BegMa:sem} in which we studied algebras that were
associative to the required order in a deformation parameter $\hbar$
but allowed the possibility that the exterior algebra in
noncommutative geometry could be nonassociative to that order. We
showed that this was {\em necessary} for the standard quantum groups
$\C_q[G]$, i.e. these associative algebras admit no associative
exterior algebra of classical dimensions that is bicovariant.
Nonassociative calculi were, however, possible by use of Drinfeld's twisiting \cite{Dri:qua} applied in the category of (super)coquasiHopf algebras. In the present work we provide many more examples
using not the quasi-Hopf algebra theory itself but a `module algebra' twist
theory in which any algebra in the category of modules covariant
under the a classical (or quantum) group is also twisted.  Such methods have
been used to obtain nonassociative algebras\cite{AlbMa:qua} as well 
as associative ones\cite{DGM:mat}. That one obtains differential calculi
as well on such algebras is explored in general terms in \cite{AkrMa:bra}. 

We show now that this setting also allows to obtain associative algebras and induced differential calculi 
for some very standard and {\em not}-quantum-group-related quantizations, but with a similar price to pay. Thus, we use Hopf algebra
methods but apply them to classical situations, notably to coadjoint
spaces $\cg^*$ and their
quantisation by the enveloping algebra $U_\hbar(\cg)$. Clearly this and other `noncommutative coordinate algebras' that we consider are perfectly associative so it is  some surprise, and  the main result of the present paper, that their natural induced noncommutative differential calculus is again
nonassociative. In the case of $U_\hbar(\cg)$ we show (Theorem~5.1.2) that any calculus which is translation and $\cg$-covariant and has classical dimensions must be nonassociative. This is analogous to the result in \cite{BegMa:sem} for quantum groups but now for classical enveloping algebras (and the proof is similar). In this way we confirm and provide major new examples of the
general analysis in \cite{BegMa:sem}. We particularly analyse the
semiclassical level of these results in terms of Poisson and
symplectic geometry followed by the next-to-semiclassical order.

An outline of the paper is as follows. In Section 2 we describe the
general algebraic twisting theory that we shall use. Section 3 then describes
the special case that will used for all our examples, namely a method of quantisation induced by a classical symmetry and a cochain.  Thus we begin with a
{\em classical} manifold $M$ with a classical Lie algebra symmetry group $\CL\subseteq {\rm diff}(M)$. 
As Hopf algebra we take $H=U(\CL)$ the enveloping algebra. Then the scheme is
that any suitable element $F\in H\tens H$ (a cochain) induces a quantisation
of $M$. We semiclassicalise this theory and  
see how Poisson-compatible (pre)connections in the sense of \cite{BegMa:sem} arise out
of the choice of $F$ and $\CL$.   The
choice of the latter covariance Lie algebra determines what kind of
connections or preconnections can arise by the cochain twisting
construction and hence what structures the quantisation respects. We
also briefly discuss the inverse problem of obtaining
a cochain $F$ and hence a quantisation given a symplectic form and
symplectic connection on $M$. Section~3.2 analyses the situation for
$M=\R^{2n}$ with its standard symplectic structure and general
symplectic connection.

Sections 4,5,6,7 then turn to  the main examples of the paper. These examples are all constructed by a second order or in some case third order analysis, i.e. we obtain the required cochain at least up to and including $\hbar^2$ terms. This will already be a substantial amount of work and is enough to expose the main phenonema. Moreover, the existence of a cochain to all orders is not really in doubt in view of the Kontsevitch universal quantisation theorem (our cochain amounts to choosing a natural `lifting' of that); our results constitute a natural choice at low order and suggest that a natural choice should be possible  to all further orders. 

We start these examples with Section~4 in which the sphere $S^2$ has a natural
cochain $F$ for covariance Lie algebra $\CL=so(1,3)$. The action of
the Lorentz group   that we use is the one on the 'sphere
at infinity' in 4-dimensional Minkowski space.  We show that one
obtains an associative quantisation of the sphere at least to  $O(\hbar^3)$ and that this
coincides with the Fedosov quantisation to this order for the
standard  Levi-Civitia connection on the sphere (which is symplectic).

Section 5 is the main example of interest in the paper. We show that  the classical enveloping algebra $U_\hbar(\cg)$ viewed as a
quantisation of $S(\cg)=\C[\cg^*]$ (functions on $\cg^*$ with its
Kirillov-Kostant bracket) can be viewed as a module-algebra cochain
twist and that this quantizes a canonical covariant preconnection in the Poisson geometry of $\cg^*$ (we show that this is in fact the only such preconnection for  all simple $\cg$ other than $sl_n$, $n>2$ and even there it is the natural choice).  The background covariance we use is $\CL=\cg\rcross \cg^*$ and we find a suitable $F$ as a powerseries  to $O(\hbar^3)$ and find that it is essentially unique to this level when we demand a further condition (Section~5.4) whereby $S(\cg^*)=\C[\cg]\subset U(\CL)$ twists into a local version of the {\em group} coordinate algebra  $\C[G]$ (see below). In effect, we require that $F$ implements the Campell-Baker-Hausdorf formula by conjugation in addition to its other properties. In Section~5.5 we discuss the Duflo map in this context and argue that the reduced form of $F$ should be the coboundary of the Duflo operator (and hence known to all orders). Although our specific universal $F$ is only found to $O(\hbar^3)$ it seems likely that these various features should extend and characterise it completely. This would be a topic for further work beyond our methods here. In Section 5.6~ we demonstrate the theory on
$\cg=\R\lcross \R=b_+$ the solvable Lie algebra in 2-dimensions. A
version of its enveloping algebra has been proposed as
 'noncommutative spacetime'\cite{MaRue:bic} and we exhibit a
  (non-unique) $F$ explicitly to  $O(\hbar^4)$ in this case.

Section 6  completes our trio of conventional examples with the Mackey
quantisation $C^\infty(N)\lcross U_\hbar(g)$ as a cochain quantisation of
$C^\infty(N)\tens S(g)\subseteq C^\infty(N\times g^*)$. This extends the model in Section~5 but we need
an extended cocycle and covariance Lie algebra $\CL=g\rcross g^*\oplus\cg$ in order to achieve this.
Section~6.3 includes the case $C^\infty(G)\lcross U_\hbar(g)$ as a
quantisation of $T^*G=G\times g^*$, where $g$ is the Lie algebra of a Lie group $G$.

We follow these with the more technical example Section~7 from quantum
group theory, which is  simply Drinfeld's theory for quantum groups $\C_q[G]$ reworked as a cochain twist. Here 
$\CL=g\oplus g^{\rm op}$ acting from the left and right and
$F,\Phi$ are built from Drinfeld's ones relating to the KZ-equations. This example is not
fundamentally new but provides the role model for our view of the more conventional
quantisations in the paper, so is included for completeness.

 Section 8 turns to the hidden nonassociativity that we have identified in the associative quantum algebras above. The most important of the many implications resulting from
the cochain twist is in Section~8.2, namely the corresponding differential calculi.
Because in our examples $F$ is not a cocycle, the exterior algebra
obtained likewise by twisting is not necessarily associative, and we show
that indeed it is not for our various examples. We describe the nonassociative differentials for each
of our examples to order  $O(\hbar^2)$. One example is a more covariant but nonassociative differential calculus for the non-commutative spacetime in \cite{MaRue:bic}. Section~8.3 shows how the same philosophy can be used to construct Dirac operators in the sense of generalised `spectral triples'. The slight generalisation beyond the axioms in \cite{Con} reflects the nonassociativity. We show that
such deformations are isospectral, a point of view consistent with other approaches such as \cite{Landi}.

 It is also true that under the cochain quantisation scheme the original
covariance becomes a (quasi)quantum group $H_F$ covariance, which we describe to $\hbar ^2$ in Section~8.1 for each of our examples. In the case of $U_\hbar(g)$ it appears by accident to be a usual (not quasi) quantum group and to be a local version of the quantum double $D(U(g))=U(g)\rcross \C[G]$, which is known to be a covariance quantum group of $U(g)$. When $U_\hbar(su_2)$ is viewed as noncommutative $\R^3$ (the so-called universal fuzzy-sphere), for example, the quantum double plays the role of quantum Euclidean group\cite{BatMa:non} motivated from 2+1 quantum gravity.  In this case the curvature of the canonical preconnection or the fact that $F$ cannot be taken to be a cocycle represents an anomaly in this quantisation of $\R^3$. The associativity obstruction in this case can in fact be resolved by adjoining an extra `time' variable and has been proposed \cite{Ma:tim} as an origin of time in noncommutative differential geometry.

Let us say finally that we work in a deformation-theoretic setting with all deformed expressions given by power-series in a parameter $\hbar$ and otherwise over $\C$;  all constructions can be formulated more (co)algebraically over any field using comodules which would, however,  be less familiar to most readers. In the main `examples' sections we work only to lower degrees in $\hbar$ for which purposes one may regard $\hbar$ as a real parameter with the deformed product  of smooth functions assumed to have these first terms in an expansion. The authors would like to thank F.W.\ Clarke for his assistance with some of the MATHEMATICA
calculations underlying the paper and Y. Bazlov for drawing our attention to the Duflo map.

\section{Preliminaries: module-algebra  cochain twists}\label{prelim}

We begin with some well known algebraic constructions, see for example the text
\cite{Ma:book}. We will only need here the classical case $H=U(\CL)$ where $\CL$ is
a Lie algebra and $A=C^\infty(M)$ where $M$ is a manifold, as the basis of the
quantisation method. In that sense we use quantum group methods but the reader does
not really need to know quantum group theory in any detail. This approach to quantisation
as been recently used in \cite{AkrMa:bra, Ma:gau} for except that
in the present paper the quantisation remains associative.

Given a Hopf algebra $(H,S,\Delta,\epsilon)$ and an invertible
$F\in H\tens H$ with $(\epsilon\tens \id)F=(\id\tens\epsilon)F=1$,
we can define a quasi-Hopf algebra $H_F=
(H,\phi,S_F,\Delta_F,\epsilon,\alpha_F,\beta_F)$, with the same
algebra and counit as
$H$, by \cite{Dri:qua}
\begin{eqnarray}
&&\Delta_F h \,=\,F.\Delta h.F^{-1}\ ,\quad
\phi\,=\, (1\tens F).(\id\tens \Delta)F.(\Delta\tens
\id)F^{-1}.(F\tens 1)^{-1}\ ,\cr
&& S_F\,=\,S\ ,\quad \alpha_F\,=\,(SF^{-(1)}).F^{-(2)}\ ,\quad
\beta_F\,=\, F^{(1)}SF^{(2)}\ .
\end{eqnarray}
In addition if there is a quasitriangular structure $\mathcal{R}$ for
$H$,
then $\mathcal{R}_F=F_{21}\mathcal{R}F^{-1}$ for $H_F$. We will call
such an $F\in H\tens H$ a 2-cochain in general, and a 2-cocycle if
$\phi=1\tens 1\tens 1$.

The significance of the twisting construction is \cite{Ma:tan} that
it corresponds to an equivalence of categories. Thus, the category
${}_{H_F}\mathcal{M}$ of left modules over $H_F$ is a monoidal
category with tensor product
operation $\tens^F$ is defined using $\Delta_F$. If
$\phi=\phi^{(1)}\tens \phi^{(2)}\tens\phi^{(3)}$, the associator
in the category  is
$\Phi_{VWZ}((v\tens w)\tens z)=\phi^{(1)}\la v\tens (\phi^{(2)}\la
w\tens \phi^{(3)}\la z)$.  On the other hand, this category is
equivalent to the category ${}_H\mathcal{M}$ of left modules over $H$
via the functor
$\mathcal{T}:{}_H\mathcal{M}\to {}_H\mathcal{M}_F$ which is
just the identity on left $H$ modules and on morhphism. A monoidal
functor also comes by definition with a natural
transformation $\vartheta:\mathcal{T}V\tens^F \mathcal{T}W\to
\mathcal{T}(V\tens W)$, given here by
$\vartheta(\mathcal{T}(v)\tens^F \mathcal{T}(w))
=\mathcal{T}(F^{-(1)}\la v\tens F^{-(2)}\la w)$ \cite{Ma:tan}.
In this way, twisting the Hopf algebra by $F$ deforms the entire
category of modules and as such deforms any and all constructions in
the category. This is the systematic 'twisting approach' to
deformation quantisation that we use.

  In particular, consider an algebra $A\in {}_H\mathcal{M}$. This
includes the
requirement that multiplication $:A\tens A\to A$ is a morphism in the
category, i.e the product is $H$-covariant (or $A$ is an $H$-module
algebra). Applying the above functor $\mathcal T$ immediately deforms
the algebra to the same vector space $A_F=A$ and the product as a map
$\mathcal{T}(A\tens A)\to \mathcal{T}A$. Using the above natural
transformation this implies a deformed product map
making an algebra $A_F\in {}_H\mathcal{M}_F$ with multiplication
$a\bullet b=(F^{-(1)}\la a)(F^{-(2)}\la b)$, and this is associative in the
category as the image of the associativity law in the undeformed
category.  This module algebra cochain quantisation method was introduced 
in \cite{DGM:mat} and related papers at the time. Examples in the cocycle case
also abound, e.g. \cite{MaOec:twi}, but the cocycle case is not 
what is of interest in the present paper since in this case the associator $\phi$ is trivial.
Neither case of `module algebra  twist' should be confused with Drinfeld's twist $H_F$ of the Hopf algebra $H$ itself. 

One may go further and  consider also the category
$\mathcal{A}\subset {}_H\mathcal{M}$ of $A$-bimodules,  which also
have $H$-actions so that the multiplications
$A\tens V\to V$ and $V\tens A\to V$
  preserve the $H$-action for all $V\in {}_A\mathcal{M}$, and so on.
Here we deform the multiplications by
$a\bullet v=(F^{-(1)}\la a).(F^{-(2)}\la v)$ and $v\bullet a=(F^{-(1)}\la
v).(F^{-(2)}\la a)$
for all $v\in V$ and $a\in A$. Similarly, if $\Omega(A)$ is an
$H$-covariant differential calculus in the sense of noncommutative
geometry (so there is for example an exterior derivative $\extd:A\to
\Omega^1(A)$ where the latter is an $A$-bimodule and $\extd$ obeys
the Leibniz rule, etc. and all
maps are morphisms in ${}_H\mathcal{M}$) then twisting any products
by the action of $F^{-1}$ gives a calculus $\Omega(A_F)$ covariant
under $H_F$. This was used for
example in \cite{AkrMa:bra}.

To this existing theory we now add some first remarks needed for the semiclassical
analysis.  As mentioned, the above should be understood as extended over formal power-series in a parameter $\hbar$ or one may continue more algebraically (using a comodule twist version of the theory). Either way, we suppose that $F^{-1}$
  is expanded as a series
\begin{eqnarray} \label{forF}
F^{-1} &=& 1\tens 1+\hbar \,  G^{(1)}+\hbar^2\,  G^{(2)}+O(\hbar^3)\ .
\end{eqnarray}
This can be inverted to give
\begin{eqnarray} \label{forFinv}
F &=& 1\tens 1-\hbar \,  G^{(1)}+\hbar^2\,
((G^{(1)})^2-G^{(2)})+O(\hbar^3)\ .
\end{eqnarray}
We can then compute
\begin{eqnarray*}
\phi &=& (1\tens F).(\id\tens \Delta)F.(\Delta\tens
\id)F^{-1}.(F\tens 1)^{-1} \cr
&=& 1\tens 1\tens 1 +\hbar\,[(\Delta\tens \id)G^{(1)}+G^{(1)}\tens
1-1\tens G^{(1)}
-(\id\tens \Delta)G^{(1)}] \cr
&&+\hbar^2\,[(\Delta\tens \id)G^{(2)}+G^{(2)}\tens 1
+1\tens ((G^{(1)})^2-G^{(2)}) +(\id\tens
\Delta)((G^{(1)})^2-G^{(2)})\cr
&&+ (\Delta\tens \id)G^{(1)}.(G^{(1)}\tens 1)
-(\id\tens \Delta)G^{(1)}.(G^{(1)}\tens 1)
-(\id\tens \Delta)G^{(1)}.(\Delta\tens \id)G^{(1)} \cr
&&-(1\tens G^{(1)}).(G^{(1)}\tens 1)
-(1\tens G^{(1)}).(\Delta\tens \id)G^{(1)}
+(1\tens G^{(1)}).(\id\tens \Delta)G^{(1)} ]         +O(\hbar^3).
\end{eqnarray*}

Now consider a special case, where $G^{(1)}=\sum X\tens Y$ (we will
later suppress the summation sign), $\Delta X=1\tens X+X\tens 1$ and
$\Delta Y=1\tens Y+Y\tens 1$. Then the order $\hbar$ part of $\phi$ is
\begin{eqnarray*}
1\tens X\tens Y+X\tens 1\tens Y+X\tens Y\tens 1 - 1\tens X\tens Y
-X\tens Y\tens 1
-X\tens 1\tens Y\,=\,0\ .
\end{eqnarray*}
The contribution from $G^{(1)}$ to the order $\hbar^2$ part of
$\phi$, using tildes to distinguish different copies of $G^{(1)}$, is
\begin{eqnarray*}
&& 1\tens X\tilde X\tens Y\tilde Y +X\tilde X\tens Y\tilde Y\tens 1
+X\tilde X\tens 1\tens Y\tilde Y+X\tilde X\tens Y\tens \tilde Y+X\tilde X\tens
\tilde Y\tens Y\cr
&&+ X\tilde X\tens \tilde Y\tens Y + \tilde X\tens X\tilde Y\tens Y
-X\tilde X\tens \tilde Y\tens Y-X\tilde X\tens Y\tilde Y\tens 1 \cr
&&-X\tilde  X \tens Y\tens \tilde Y - X\tilde  X \tens 1\tens Y \tilde Y
-X\tens Y\tilde X\tens\tilde Y- X\tens \tilde X\tens Y\tilde Y - \tilde X\tens
X\tilde Y\tens Y\cr
&&-\tilde X\tens X\tens Y\tilde Y - 1\tens X\tilde X\tens Y \tilde Y
+ \tilde X\tens X\tilde Y \tens Y + \tilde X\tens X\tens Y\tilde Y\cr
&=& X\tilde X\tens \tilde Y\tens Y+ \tilde X\tens X\tilde Y \tens Y
-X\tens Y\tilde X\tens\tilde Y- X\tens \tilde X\tens Y\tilde Y \cr
&=& X\tilde X\tens \tilde Y\tens Y+ \tilde X\tens X\tilde Y \tens Y
-\tilde X\tens \tilde Y X\tens Y- X\tens \tilde X\tens Y\tilde Y \cr
&=& X\tilde X\tens \tilde Y\tens Y+ \tilde X\tens [X,\tilde Y] \tens Y
- X\tens \tilde X\tens Y\tilde Y\ .
\end{eqnarray*}
To summarise, if we put, for $P\in H\tens H$,
\begin{eqnarray*}
\del P &=& 1\tens P -(\Delta\tens\id)P+(\id\tens\Delta)P-P\tens 1 \ .
\end{eqnarray*}
then the expression for $\phi$ becomes
\begin{eqnarray*}
\phi \,=\, 1\tens 1\tens 1+\hbar^2(X\tilde X\tens \tilde Y\tens Y+ \tilde
X\tens [X,\tilde Y] \tens Y
- X\tens \tilde X\tens Y\tilde Y-\del G^{(2)})+O(\hbar^3)\ .
\end{eqnarray*}
Next note that, if $A,B,C,D\in H$ are also primitive (have linear coproducts
like $X,Y$) then  
\begin{eqnarray*}
\del(AB\tens CD) &=& AB\tens C\tens D+AB\tens D\tens C -
A\tens B\tens CD-B\tens A\tens CD\ .
\end{eqnarray*}
Hence if we were to put $G^{(2)}=(G^{(1)})^2/2$, then
\begin{eqnarray}\label{iuyf}
\phi \,=\, 1\tens 1\tens 1+\hbar^2([X,\tilde X]\tens \tilde Y\tens Y+
2\,\tilde X\tens [X,\tilde Y] \tens Y
- X\tens \tilde X\tens [Y,\tilde Y])/2+O(\hbar^3)\ .
\end{eqnarray}

  The set of primitive elements in $H$ form a Lie algebra in $H$ which we can view as a Lie
bialgebra with zero Lie cobracket. After twisting these elements
acquire a Lie cobracket
  \[ \delta Z=\sum [X,Z]\tens Y+X\tens [Z,Y]-[Y,Z]\tens X-Y\tens
[X,Z]\]
  for all $Z$ in the Lie algebra and together with the $\hbar^2$ part
of $\phi$ form a quasi-Lie bialgebra. This is the infinitesimal object associated to the quasi-quantum group $H_F$. 

\section{Poisson-compatible connections from cochain twists at the
semiclassical level}
\label{bvcihjvx}

We start by briefly recalling the main ideas of \cite{BegMa:sem}. As is well
known, if one considers the quantisation of the functions
$C^\infty(M)$ of a classical manifold $M$, the initial data usually
specified is a Poisson structure defined by a bivector $\omega$ (in
the symplectic case this is invertible with inverse (also denoted
$\omega$) a closed 2-form). Any flat deformation-quantisation
$A_\hbar$ corresponds on looking at the leading part of its
commutator
\eqn{commutator}{ a\bullet b-b\bullet a=\hbar \{a,b\}+ O(\hbar^2)}
to a Poisson bracket $\{a,b\}=\omega(\extd a,\extd b)$. More
recently, we considered the same question for a
noncommutative differential calculus $\Omega(A_\hbar)$ quantizing the usual
exterior algebra $\Omega(M)$ but in a slightly
weaker than usual setting (without assuming associativity of products
involving differential forms). We found that the initial data for
this at least in the symplectic case was a compatible connection $\nabla$ defined by
\eqn{nablahata}{a\bullet \extd b-\extd b\bullet a=\hbar\nabla_{\hat a}\extd
b+O(\hbar^2).}
Here $\hat a$ denotes the Hamiltonian vector field $\omega(\extd a,\
)=\{a,\ \}$. Here $\nabla$ is not necessarily well-defined in the Poisson case  even along Hamiltonian vector fields; it could be called a partial connection where defined or we should speak more precisely of a 'preconnection' $\varnabla$ defined almost identically by
 \eqn{prenabla}{a\bullet \extd b-\extd b\bullet a=\hbar\varnabla_{a}\extd
b+O(\hbar^2)}
(here $\varnabla$ was called $\gamma$ in \cite{BegMa:sem}).  (In fact there is a more general notion of `contravariant connection' which can also be used here, see\cite{Haw}).  The (Poisson)-compatibility condition is
\eqn{compat}{\nabla_{\hat a}\extd b-\nabla_{\hat b}\extd a=\extd\{a,b\}}
and under some mild conditions in the symplectic case becomes \cite{BegMa:sem} that
$\nabla$ is a torsion free symplectic connection in the usual
sense. Finally, the curvature and torsion of the connection 
\eqn{RT}{R_\nabla(\hat a,\hat b)=[\nabla_{\hat a},\nabla_{\hat b}]-\nabla_{[\hat a,\hat b]},\quad T_\nabla(\hat a,\hat b)=\nabla_{\hat a}\hat b-\nabla_{\hat b}\hat a-[\hat a,\hat b]}
are defined in the usual way as for any conneciton. In terms of a preconnection the equations are more precisely
\eqn{preRT}{ \varnabla_{a}\extd b-\varnabla_{b}\extd a=\extd\{a,b\},\quad R_{\varnabla}(a,b)=[\varnabla_a,\varnabla_b]-\varnabla_{\{a,b\}},\quad T_\varnabla(a,b)=\varnabla_a\hat b-\varnabla_b\hat a-[\hat a,\hat b]}
where the last term in the curvature is in view of $[\hat a,\hat b]=\widehat{\{a,b\}}$.  It was shown in \cite{BegMa:sem} that the curvature coincides with the Jacobiator or obstruction to associativity for the
differential calculus $\Omega(A_\hbar)$ at the relevant lowest order. This was also observed 
recently in  \cite{Haw} where it was shown that the associative case corresponds to a zero-curvature (contravariant) connection, although we were not aware of this at the time of  \cite{BegMa:sem}. 

{}From a geometrical point of view, however, it would seem at the
semiclassical level quite reasonable to consider symplectic or Poisson manifolds
equipped with  connections with curvature, which in
quantisation terms would mean by \cite{BegMa:sem} associative quantum
algebras with nonassociative differential calculi. We have seen in
Section~2 a general method to construct examples of such hybrid
situations by means of cochain twists. Now we see what that amounts
to at the semiclassical level.

\subsection{Inducing connections by twisting}

We consider the diffeomorphism group acting on the functions on a
manifold $M$.
This action extends to the vector fields and forms on the manifold,
and infinitesimally the action is called the Lie derivative.
A vector field $X$ (i.e.\ an element of the Lie algebra of the
diffeomorphism group)
acts on forms by $\mathcal{L}_X\xi=\varpi_X(\extd
\xi)+\extd(\varpi_X\xi)$,
where $\varpi_X$ is the interior product.  If $\xi=\xi_{k_1\dots
k_n}\extd x^{k_1}\wedge\dots\wedge\extd x^{k_n}$, then
\begin{eqnarray*}
\mathcal{L}_X\xi &=& X^p\,\xi_{k_1\dots k_n,p}\,
\extd x^{k_1}\wedge\dots\wedge\extd x^{k_n}
+\xi_{k_1\dots k_n}\extd(\varpi_X(\extd x^{k_1}\wedge\dots\wedge\extd
x^{k_n})) \cr
&=&  \Big(X^j\,\xi_{k_1\dots k_n,j}\,
  +X^j_{\phantom{j},k_i}\,
\xi_{k_1\dots k_{i-1},j,k_{i+1}\dots k_n}\Big)
\extd x^{k_1}\wedge\dots\wedge\extd x^{k_n}
\end{eqnarray*}
We , can therefore apply the ideas of Section~2 with $A=C^\infty(M)$
and $H=U(\diff(M))$ acting via the Lie derivative on all tensorial
objects. We will proceed  with $F,F^{-1}$ power-series having
values in $U({\rm diff}(M))\tens U({\rm diff}(M))$, however we are interested 
in this section only in the differential geometry resulting from the semiclassical part and not
in the formal construction of these objects. 

We assume an expansion of the 2-cochain
$F^{-1}=\id^2+\hbar\,G^{(1)}+O(\hbar^2)$, where
  $G^{(1)}=\sum X\tens Y\in {\rm diff}(M)\tens {\rm diff}(M)$, acts
on $\Omega^n(M)\tens \Omega^m(M)$ by $\xi\tens\eta\mapsto
\xi\tens\eta+\hbar\,
\sum \mathcal{L}_X\xi\tens \mathcal{L}_Y\eta+O(\hbar^2)$. Then
the commutator of two functions $a,b\in C^\infty(M)$ is
\begin{eqnarray} \label{kjcvjfghxx}
\{a,b\} &=& \hbar\sum\Big((\mathcal{L}_X a)(\mathcal{L}_Y b)-
(\mathcal{L}_Y a)(\mathcal{L}_X b)\Big)+O(\hbar^2) \cr
&=& \hbar\,\sum (X^i\,Y^j-X^j\,Y^i)\,a_{,i}\,b_{,j}+O(\hbar^2)\ .
\end{eqnarray}
Now we set $\omega=\sum( X\tens_M Y- Y\tens_M X)\in
{\rm diff}(M)\tens_M {\rm diff}(M)$ (putting
summation subscripts on $X$ and $Y$
would only be confusing). This is antisymmetric, so we have
$\omega\in {\rm diff}(M)\wedge {\rm diff}(M)$. We shall assume for
convenience that
we are in the symplectic case, where $\omega$ is invertible, and its
inverse is
a closed 2-form, otherwise more generally we assume that $\omega$ is
a Poisson bivector, i.e. induces a Poisson bracket. A sufficient
condition for the latter is that $\phi=1\tens1\tens 1$ when projected
from $\tens$ to  $\tens_M$. Another sufficient condition is that
$G^{(1)}$ obeys the Classical Yang-Baxter equations, but neither is
necessary and we do not assume them.

Now apply the same $F^{-1}$ to deform the products of functions and
1-forms as explained in Section~2. This implies a connection $\nabla$
resulting from
the commutator of a function and a 1-form:
\begin{eqnarray*}
[a,\xi_i\,\extd x^i] &=& \hbar\sum X^k\, a_{,k}(\xi_{i,j}\,Y^j\,\extd
x^i+\xi_i\,Y^i_{\phantom{i},j}\,
\extd x^j)-(X\leftrightarrow Y)+O(\hbar^2) \cr
&=&\hbar \sum (X^k\,Y^j-X^j\,Y^k)\,a_{,k}\,\xi_{i,j}\,\extd x^i+
\sum(X^k\,Y^i_{\phantom{i},j}-Y^k\,X^i_{\phantom{i},j})\,\xi_i\,a_{,k}\,\extd
x^j +O(\hbar^2)\cr
&=& \hbar\,\omega^{kj}\,a_{,k}\,\xi_{i,j}\,\extd x^i+
\sum(X^k\,Y^i_{\phantom{i},j}-Y^k\,X^i_{\phantom{i},j})\,\xi_i\,a_{,k}\,\extd
x^j +O(\hbar^2)\cr
&=& \hbar\,\omega^{kj}\,a_{,k}\Big(\xi_{i,j}\,\extd x^i+
\sum \omega_{js}\,
(X^s\,Y^i_{\phantom{i},p}-Y^s\,X^i_{\phantom{i},p})\,\xi_i\,\extd
x^p\Big)+O(\hbar^2) \cr
&=& \hbar\,\omega^{kj}\,a_{,k}\Big(\xi_{i,j}\,\extd x^i-
\Gamma^i_{jp}\,\xi_i\,\extd x^p\Big)+O(\hbar^2) \ .
\end{eqnarray*}
Then the Christoffel symbols of the connection can be seen to be
\begin{eqnarray} \label{smmbolsvcwu}
\Gamma^i_{jp} &=& -\sum \omega_{js}\,
(X^s\,Y^i_{\phantom{i},p}-Y^s\,X^i_{\phantom{i},p})\ .
\end{eqnarray}
\begin{propos}
The connection is characterised by the equation, for
$a\in C^\infty(M)$ and $\xi\in\Omega^1 M$:
\[ \nabla_{\hat a}\xi=\sum X(a)\,(\extd \<Y,\xi\>+\varpi_Y\extd\xi)
-Y(a)\,(\extd \<X,\xi\>+\varpi_X\extd\xi)\ .\]
\end{propos}
\noindent {\bf Proof:}\quad Just substitute from (\ref{smmbolsvcwu})
for the Christoffel
symbols.\eproof

\medskip

The condition that $\omega$ is closed is
\begin{eqnarray*}
\omega_{ir}\,\omega^{rs}_{\phantom{rs},k}\,\omega_{sj}\,\extd
x^k\wedge
\extd x^i\wedge \extd x^j\,=\,0\ .
\end{eqnarray*}
From \cite{BegMa:sem} the connection is necessarily comapatible with the
differential
structure in the form
\begin{eqnarray}
\frac{\partial \omega^{ij}}{\partial x^p} &=&
\omega^{jq}\,\Gamma^i_{qp}-
\omega^{iq}\,\Gamma^j_{qp}\ .
\end{eqnarray}
For the Christoffel symbols in (\ref{smmbolsvcwu}) we have
\begin{eqnarray*}
\omega^{jq}\,\Gamma^i_{qp}-
\omega^{iq}\,\Gamma^j_{qp} &=&
\sum
(X^i\,Y^j_{\phantom{i},p}-Y^i\,X^j_{\phantom{i},p})
- \sum
(X^j\,Y^i_{\phantom{i},p}-Y^j\,X^i_{\phantom{i},p}) \cr
&=& \sum(X^i\,Y^j-X^j\,Y^i)_{,p}\,=\,\omega^{ij}_{\phantom{ij},p}\ ,
\end{eqnarray*}
as expected.  

In summary, any manifold may potentially by quantised by choosing a
cochain $F$ with values in $U(\diff{M})\tens U(\diff{M})[[\hbar]]$. The
leading order of $F$ which we have denoted $\sum X\tens Y$ will
induce a bivector which will not necessarily be a Poisson bivector
(the quantised algebra may not necessarily be associative). However,
if it is, we will also have induced a Poisson-compatible connection (or more 
precisely a preconnection). Conversely, given, say, a symplectic manifold $M$ equipped with
 (as in Fedosov theory) a symplectic
connection $\nabla$ we can look for a suitable $F$ inducing these
intial data to lowest order and such that $\phi=1\tens 1\tens 1$ when
tensorised over $\tens_M$. This provides an alternative and more
categorical approach to the quantisation problem in the spirit of
Fedosov theory but now having the merit of also quantising
differential calculi and all other covariant constructions, albeit
with potential nonassociativity.

\subsection{The inverse problem for $\R^{2n}$}

The `inverse problem' of finding $F$ even to lowest order (i.e. $\sum
X\tens Y$) such that a given Poisson bi-vector $\omega$ is obtained
and a given symplectic or Poisson-compatible (pre)connection is obtained
as above appears to be a tricky one. Here we will look at what is
involved in the simplest possible case of $\R^{2n}$. We take the
standard symplectic structure and note that torsion free
  symplectic connections for it are in 1-1 correspondence with totally
symmetric Christoffel symbols
$\Gamma_{abc}=\omega_{ad}\,\Gamma^d_{bc}$ \cite{GRS:Fed}. 

  We can easily make the canonical symplectic form for $\R^{2n}$ by
adding
  $\hbar X\tens Y$ terms to $F^{-1}$ where $X,Y$ are constant vector
fields.
  By (\ref{smmbolsvcwu}) these will give zero Christoffel symbols. But
then we can add
  further terms to $F^{-1}$ of the form $\sum X\tens Y =\sum f.U\tens
V-U\tens f.V$, where
  $f$ is a real valued function and $U,V$ are constant vectors. Then
  \begin{eqnarray*}
\Gamma_{ijp} &=& -\sum \omega_{ik}\,\omega_{js}\,f_{,p}\,
(-U^s\,V^k-V^s\,U^k)\ .
\end{eqnarray*}
If we set $f$ to be the linear function $\omega_{pq}\,W^q\,x^p$, then
  \begin{eqnarray*}
\Gamma_{ijp} &=& \omega_{ik}\,\omega_{js}\,\omega_{pq}\sum \,W^q\,
(U^s\,V^k+V^s\,U^k)\ .
\end{eqnarray*}
By adding terms of this form we can recreate any symplectic
connection with
constant Christoffel symbols by this form of $F$. 
 Then the curvature is given by
\begin{eqnarray*}
R^a_{\phantom{a}bcd} &=& \omega^{me}\,\omega^{ag}(\Gamma_{edb}\Gamma_{gcm}
-\Gamma_{ecb}\Gamma_{gdm}) \cr 
&=& \omega^{me}\,\omega^{ag}\,
\omega_{ek}\,\omega_{bq}\,
\omega_{gp}\,\omega_{mr} \Big(\omega_{ds}\,
\omega_{cv}\sum \,W^q\,
(U^s\,V^k+V^s\,U^k) \,\tilde W^r\,
(\tilde U^v\,\tilde V^p+\tilde V^v\,\tilde U^p) \cr
&&-\, \omega_{cs}\,\omega_{dv}\sum \,W^q\,
(U^s\,V^k+V^s\,U^k) \,\tilde W^r\,
(\tilde U^v\,\tilde V^p+\tilde V^v\,\tilde U^p)\Big) \cr
&=&  \omega_{bq}\,\omega_{kr} (\omega_{ds}\,
\omega_{cv}-\, \omega_{cs}\,\omega_{dv})\sum \,W^q\,
(U^s\,V^k+V^s\,U^k) \,\tilde W^r\,
(\tilde U^v\,\tilde V^a+\tilde V^v\,\tilde U^a)  \cr
&=&  \omega_{bq}\,\omega_{kr} (\omega_{ds}\,
\omega_{cv}-\, \omega_{cs}\,\omega_{dv})\sum \,W^q\,
(U^s\,V^k\tilde V^v\,\tilde U^a+V^s\,U^k\tilde U^v\,\tilde V^a) \,\tilde W^r \ ,
\end{eqnarray*}
where the tildes denote a second set of triples $(U,V,W)$ and the sum is over both sets.

\section{Quantising $S^2$ by cochain twist}

Here we describe a simple example of the cochain quantisation method in
Section \ref{bvcihjvx}.
The covariance used for the twisting will be the Lorentz group and
its action on $S^2$, a  nonlinear one related to spacetime physics (the sphere at infinity in
Minkowski space).  This induces a quantisation not related as far as we know to
the representation theoretic coadjoint orbit  examples given later.

\subsection{Some nice vector fields on the 2-sphere}

Our goal is to show how a natural covariance, cochain and hence quantisation arise 
in a nice way from the geometry of  $S^2=\{(x,y,z)\in\mathbb{R}^3:x^2+y^2+z^2=1\}$. 
We use the standard inner product on $\mathbb{R}^3$. 

Thus, given
$\underline v\in \mathbb{R}^3$, we have the natural vector field $X[\underline v]$
defined by
$X[\underline v](\underline r)
=\underline v-\underline r\,\<\underline v,\underline r\>$ (for
$\underline r\in
S^2$) which tangent to the sphere at $\underline r$. Also at each such point we have
the orbital angular moment   vector
field
$Y[\underline v](\underline r)=\underline v\times\underline r$, where
$\times$ is the vector cross product. These vector fields are clearly well behaved under rotation; consider an orthogonal
transformation $T\in O_3(\mathbb{R})$. Then $T(X[\underline
v])(\underline r)$
is by definition
\begin{eqnarray*}
T(X[\underline v](T^{-1}\underline r))
\,=\, T\underline v-T(T^{-1}\underline r)\,\<\underline
v,T^{-1}\underline r\>
\,=\, T\underline v-\underline r\,\<T\underline v,\underline r\>\,=\,
X[T\underline v](\underline r)\ .
\end{eqnarray*}
Also we have $T(Y[\underline v])(\underline r)$ equal to
\begin{eqnarray*}
T(Y[\underline v](T^{-1}\underline r))
\,=\, T(\underline v\times T^{-1}\underline r)\,=\, {\rm
det}(T).T\underline v\times
\underline r\,=\, {\rm det}(T).Y[T\underline v](\underline r)\ ,
\end{eqnarray*}
where the determinant enters by the change in sign of the vector
product under a change in orientation of ${\mathbb R}^3$.

The Lie bracket of two vector fields is defined as usual by
$[X,Y]^i=Y^i_{\phantom{i},j}\, X^j- X^i_{\phantom{i},j}\, Y^j$. Thus
\begin{eqnarray*}
[Y[\underline v],Y[\underline w]](\underline r) &=&
\underline  w \times (\underline  v\times \underline r)-
\underline v\times(\underline  w\times \underline r)  \cr
&=& \<\underline w,\underline r\>\,\underline v-
  \<\underline w,\underline v\>\,\underline r
  - \<\underline v,\underline r\>\,\underline w
  + \<\underline v,\underline w\>\,\underline r\cr
&=&  \<\underline w,\underline r\>\,\underline v
  - \<\underline v,\underline r\>\,\underline w\ ,\cr
[Y[\underline v],X[\underline w]](\underline r) &=&
-\underline r\,\<\underline w,Y[\underline v](\underline r)\>
-Y[\underline v](\underline r)\,\<\underline w,\underline r\>
-\underline v\times X[\underline w](\underline r) \cr
&=&
-\underline r\,\<\underline w,\underline  v\times \underline r\>
-\underline  v\times \underline r\,\<\underline w,\underline r\>
-\underline v\times (\underline w-\underline r\,\<\underline
w,\underline r\>) \cr
&=&
-\underline r\,\<\underline w,\underline  v\times \underline r\>
-\underline v\times \underline w \,=\,
\underline r\,\<\underline  v\times \underline w,\underline r\>
-\underline v\times \underline w\cr
&=& - X[\underline v\times \underline w ](\underline r) \ ,\cr
[X[\underline v],X[\underline w]](\underline r) &=&
-\underline r\,\<\underline w,X[\underline v](\underline r)\>
-X[\underline v](\underline r)\,\<\underline w,\underline r\>
+\underline r\,\<\underline v,X[\underline w](\underline r)\>
+X[\underline w](\underline r)\,\<\underline v,\underline r\> \cr
&=&
-\underline r\,\<\underline w,\underline v\>
-\underline v\,\<\underline w,\underline r\>
+\underline r\,\<\underline v,\underline w\>
+\underline w\,\<\underline v,\underline r\> \cr
&&
+\underline r\,\<\underline w,\underline r\,\<\underline v,\underline
r\>\>
+\underline r\,\<\underline v,\underline r\>\,\<\underline
w,\underline r\>
-\underline r\,\<\underline v,\underline r\,\<\underline w,\underline
r\>\>
-\underline r\,\<\underline w,\underline r\>\,\<\underline
v,\underline r\> \cr
&=&
-\underline v\,\<\underline w,\underline r\>
+\underline w\,\<\underline v,\underline r\> 
\end{eqnarray*}
whereas
\begin{eqnarray*}
Y[\underline v\times \underline w](\underline r)
  &=& (\underline v\times \underline w)\times \underline r \,=\,
  \<\underline r ,\underline v\>\, \underline w-
  \<\underline r ,\underline w\>\, \underline v\ .
\end{eqnarray*}
Hence we have
  \begin{eqnarray}
[Y[\underline v],Y[\underline w]] \,=\, -Y[\underline
v\times\underline w]\ ,\quad
[Y[\underline v],X[\underline w]]\,=\,- X[\underline v\times
\underline w ]\ ,\quad
[X[\underline v],X[\underline w]]\,=\, Y[\underline v\times\underline
w]\ .
\end{eqnarray}

In other words, the $Y$ fields generate rotations and the $X$
generare boosts of the Lie algebra $\CL=so(1,3)\subset\diff(S^2)$. This action has the physical interpretation
mentioned above and will be used to induce the quantisaton.

\subsection{The first order rotation invariant  2-cochain} Set
\begin{eqnarray*}
&& X_1\,=\, X[(1,0,0)]\ ,\quad X_2,=\, X[(0,1,0)]\ ,\quad X_3\,=\,
X[(0,0,1)]\ ,\cr
&& Y_1\,=\, Y[(1,0,0)]\ ,\quad Y_2\,=\, Y[(0,1,0)]\ ,\quad Y_3\,=\,
Y[(0,0,1)]\ .
\end{eqnarray*}
Using the matrix coefficients of $T\in O_3(\mathbb{R})$ in the
standard basis,
\begin{eqnarray*}
T(X_i)\,=\, \sum_j T_{ji}\,X_j\ ,\quad T(Y_i)\,=\, \det T.\sum_j
T_{ji}\,Y_j\ ,
\end{eqnarray*}
so under a rotation each $X_i$ is sent to a linear combination
of the $X_j$ with (the important bit) constant coefficients, not
general functions on $S^2$, and likewise with the $Y_i$.
Using these notations, for any $3\times 3$ real matrix $\kappa$,
emphasising the fact that the tensor product is over $\mathbb{R}$, we take
the lowest order part of $F,F^{-1}$ in the form
\begin{eqnarray*}
G^{(1)}_\kappa\,=\, \kappa_{ij} \, X_i\tens_\mathbb{R} Y_j .
\end{eqnarray*}
Then,
\begin{eqnarray*}
(T\tens T)G^{(1)}_\kappa \,=\, \det T.\kappa_{ij} \,T_{ki}\,
X_k\tens_\mathbb{R}
  T_{sj}\,Y_s\,=\,
\det T.T_{ki}\, \kappa_{ij}\,T^T_{js} \,X_k\tens_\mathbb{R} Y_s\, =\, \det T.G^{(1)}_{T\kappa T^{-1}}
\end{eqnarray*}
Hence if we want $\Upsilon[\kappa]$ to be rotation invariant, then we would
like
$T\kappa T^{-1}=\kappa$ for all $T\in SO_3(\mathbb{R})$; we therefore take 
$\kappa$ to be (half) the identity matrix. Of course it will still have
its sign changed by
orientation reversing $T\in O_3(\mathbb{R})$, but this is also
true of a rotation invariant symplectic form on $S^2$, so that this is exactly what
we want. This last consideration also excludes terms of the form $X\tens X$ and $Y\tens Y$ in
our ansatz for $G^{(1)}$.

We are therefore led by rotational considerations to $G^{(1)}={1\over 2} X_i\tens Y_i$, which we use henceforth. We use  (\ref{kjcvjfghxx}) to find
the corresponding
Poisson structure. As it is rotation invariant, we only have to
evaluate $G^{(1)}$ at the point $(0,0,1)$:
\begin{eqnarray*}
2 G^{(1)}(0,0,1) &=& X_1(0,0,1)\tens Y_1(0,0,1)+X_2(0,0,1)\tens
Y_2(0,0,1) \cr
&=& -(1,0,0) \tens (0,1,0) +(0,1,0) \tens (1,0,0)\ .
\end{eqnarray*}
This means that, were we to reduce to $\tens_{C(S^2)}$, we would get
the
Poisson structure corresponding to the usual symplectic form on
$S^2$.

\subsection{The connection} \label{fycy}
It will be convenient to choose coordinates $(x,y)\in\mathbb{R}^2$
for the  hemisphere $\{(x,y,z)\in\mathbb{R}^3:x^2+y^2+z^2=1
\ {\rm and}\ \ z>0\}$. Then the component vector fields are given by
\begin{eqnarray*}
&& X_1(x,y) \,=\,(1-x^2,-xy)\ ,\quad X_2(x,y) \,=\,(-xy,1-y^2)
\ ,\quad X_3(x,y) \,=\,(-xz,-yz)\ ,\cr
&& Y_1(x,y) \,=\,(0,-z)\ ,\quad Y_2(x,y) \,=\,(z,0)\ ,\quad Y_3(x,y)
\,=\,(-y,x)\ .
\end{eqnarray*}
The Poisson tensor corresponding to the unique $G^{(1)}$ found above is 
 \begin{eqnarray*}
2\,\omega &=& \sum_i \Big(X_i\tens_{C(S^2)}Y_i
-Y_i\tens_{C(S^2)}X_i\Big)\ ,
\end{eqnarray*}
so, numbering the coordiantes $x^1=x$ and $x^2=y$, $\omega^{ij}$
is antisymmetric and $\omega^{12}(x,y)=-z$. Taking the inverse matrix
gives $\omega_{12}=1/z$. From (\ref{smmbolsvcwu}) the Christoffel
symbols are
(with summation sign supressed)
\begin{eqnarray*}
2\,\Gamma^i_{1p} &=& \omega_{1s}\,
(Y_k^s\,X^i_{k,p}-X_k^s\,Y^i_{k,p}) \,=\,
(Y_k^2\,X^i_{k,p}-X_k^2\,Y^i_{k,p})/z\ ,  \cr
2\,\Gamma^i_{2p} &=& \omega_{2s}\,
(Y_k^s\,X^i_{k,p}-X_k^s\,Y^i_{k,p}) \,=\,
-(Y_k^1\,X^i_{k,p}-X_k^1\,Y^i_{k,p})/z\ .
\end{eqnarray*}
In more detail,
\begin{eqnarray*}
2\,\Gamma^1_{1p} &=&
(Y_k^2\,X^1_{k,p}-X_k^2\,Y^1_{k,p})/z \cr
&=& (-z(1-x^2)_{,p}+x(-xz)_{,p}-(1-y^2)\,z_{,p}-yz\,y_{,p})/z\ ,\cr
  2\,\Gamma^2_{1p} &=&
(Y_k^2\,X^2_{k,p}-X_k^2\,Y^2_{k,p})/z \cr
&=& (-z(-xy)_{,p}+x(-yz)_{,p}+xy(-z)_{,p}+yz\,x_{,p})/z\ ,\cr
2\,\Gamma^1_{2p} &=&
(X_k^1\,Y^1_{k,p}-Y_k^1\,X^1_{k,p})/z \cr
&=&(-xy\,z_{,p}+xz\,y_{,p}-z(-xy)_{,p}+y(-xz)_{,p})/z\ , \cr
2\,\Gamma^2_{2p} &=&
(X_k^1\,Y^2_{k,p}-Y_k^1\,X^2_{k,p})/z \cr
&=& ((1-x^2)(-z)_{,p}-xz\,x_{,p}-z(1-y^2)_{,p}+y(-yz)_{,p})/z\ .
\end{eqnarray*}
Then we get, using $z_{,1}=-x/z$ and $z_{,2}=-y/z$,
\begin{eqnarray*}
2\,\Gamma^1_{1p}
&=& (zx\,x_{,p}-x^2\,z_{,p}-(1-y^2)\,z_{,p}-yz\,y_{,p})/z\ ,\cr
2\,\Gamma^1_{11}
&=& (zx+x^3/z+(1-y^2)x/z)/z \,=\,x(z^2+x^2+1-y^2)/z^2\ ,\cr
\Gamma^1_{11} &=& x(1-y^2)/z^2\ ,\cr
2\,\Gamma^1_{12}
&=& (x^2\,y/z+(1-y^2)\,y/z-yz)/z \,=\, y (x^2+1-y^2-z^2)/z^2\ , \cr
\Gamma^1_{12} &=& x^2\,y/z^2\ ,\cr
  2\,\Gamma^2_{1p}
&=& (z(xy)_{,p}-x(yz)_{,p}-xy\,z_{,p}+yz\,x_{,p})/z
\,=\, 2(zx_{,p}y-xyz_{,p})/z\ ,\cr
\Gamma^2_{11} &=&y (1-y^2)/z^2\ , \cr
\Gamma^2_{12} &=& xy^2/z^2\ , \cr
2\,\Gamma^1_{2p}
&=&(-xy\,z_{,p}+xz\,y_{,p}+z(xy)_{,p}-y(xz)_{,p})/z \,=\, 2
(xz\,y_{,p}-yxz_{,p})/z\ , \cr
\Gamma^1_{21} &=& x^2\,y/z^2\ , \cr
\Gamma^1_{22} &=& x(1-x^2)/z^2\ , \cr
2\,\Gamma^2_{2p}
&=& (-(1-x^2)z_{,p}-xz\,x_{,p}+zyy_{,p}-y^2\,z_{,p})/z\ , \cr
2\,\Gamma^2_{22}
&=& ((1-x^2)y/z+zy+y^3/z)/z\,=\, y(1-x^2+z^2+y^2)/z^2\ , \cr
\Gamma^2_{22}
&=&   y(1-x^2)/z^2\ ,   \cr
2\,\Gamma^2_{21}
&=& ((1-x^2)x/z-xz+y^2x/z)/z \,=\, x (1-x^2-z^2+y^2)/z^2\ , \cr
\Gamma^2_{21}  &=& x\,y^2/z^2\ .
\end{eqnarray*}
  From this we see that the connection is torsion free, and since it
is also compatible with
the differential structure, by \cite[Sec. 3]{BegMa:sem} it is also symplectic.

\subsection{Metric compatability}
The metric induced from the standard embedding in $\mathbb{R}^3$
with the standard inner product is, in $(x,y)$ coordinates,
\begin{eqnarray*}
g=\{g_{ij} \}&=& \frac{1}{z^2}
\left(\begin{array}{cc}1-y^2 & xy \\xy & 1-x^2\end{array}\right)\ .
\end{eqnarray*}
If we organise the Christoffel symbols into the matrices
\begin{eqnarray*}
N_k &=& \left(\begin{array}{cc}\Gamma^1_{k1} & \Gamma^2_{k1}
  \\  \Gamma^1_{k2} &  \Gamma^2_{k2} \end{array}\right)\ ,
\end{eqnarray*}
then, using matrix multiplication,
\begin{eqnarray*}
\nabla_k g_{ij} &=& \frac{\partial g_{ij}
}{\partial x^k}
- N_k.g_{ij}    -     g_{ij}. N^T_k\ .
\end{eqnarray*}
In our case
\begin{eqnarray*}
N_1 \,=\, \frac{1}{z^2}\left(\begin{array}{cc} x(1-y^2) & y (1-y^2)
  \\  x^2\,y &  x\, y^2 \end{array}\right)\ ,\quad
  N_2 \,=\, \frac{1}{z^2}\left(\begin{array}{cc} x^2\,y & x\,y^2
  \\  x(1-x^2) &  y(1-x^2)\end{array}\right)\ ,
\end{eqnarray*}
and using this it can quickly be verified that the covariant
derivatives
of the metric vanish. Since the induced connection above was torsion free,
it must be the usual Levi-Civita connection on $S^2$.

\subsection{The curvature}
In a coordinate frame, the curvature is given by
\[
R^{l}_{\phantom{l}ijk}\,=\, \frac{\partial \Gamma^l_{ki}}{\partial
x^j}\,-\, \frac{\partial \Gamma^l_{ji}}{\partial x^k}\,+\,
\Gamma^m_{ki}\,\Gamma^l_{jm}\,-\,\Gamma^m_{ji}\,\Gamma^l_{km}\ .
\]
At the point $x=y=0$ we find that, to first order in $x$ and $y$,
$\Gamma^1_{11}=\Gamma^1_{22} = x$, $\Gamma^2_{11}=\Gamma^2_{22} = y $
and all other Christoffel symbols vanish. Then, at that point,
\[
R^{l}_{\phantom{l}ijk}\,=\, \frac{\partial \Gamma^l_{ki}}{\partial
x^j}\,-\, \frac{\partial \Gamma^l_{ji}}{\partial x^k} \,=\,
\delta^l_j\,\delta_{ki}-\delta^l_k\,\delta_{ji}\ .
\]

Now we set $G^{(1)}=X_i\tens Y_i/2$ (summed over $i$) as above, and 
$G^{(2)}=(G^{(1)})^2/2$, then from (\ref{iuyf})
\begin{eqnarray} \label{vcbbvcdiblv}
\phi &=& 1\tens 1\tens 1+\hbar^2([X_i, X_j]\tens  Y_j\tens Y_i+
2\, X_j\tens [X_i, Y_j] \tens Y_i \cr
&& - X_i\tens  X_j\tens [Y_i, Y_j])/8+O(\hbar^3)\cr
&=&  1\tens 1\tens 1+\hbar^2(Y[\underline e_i\times\underline
e_j]\tens  Y_j\tens Y_i-
2\, X_j\tens X[\underline e_i\times \underline e_j]\tens Y_i \cr
&& + X_i\tens  X_j\tens Y[\underline e_i\times \underline
e_j])/8+O(\hbar^3)\ ,
\end{eqnarray}
where the $\underline e_i$ are the usual basis vectors.
Now we can expand the summations:
\begin{eqnarray*}
X_i\tens  X_j\tens Y[\underline e_i\times \underline e_j] &=&
X_{[1]}\tens  X_{[2]}\tens Y_{[3]}-X_{[1]}\tens  X_{[3]}\tens Y_{[2]}
+X_{[2]}\tens  X_{[3]}\tens Y_{[1]} \cr
&& - X_{[2]}\tens  X_{[1]}\tens Y_{[3]}+X_{[3]}\tens  X_{[1]}\tens
Y_{[2]}
-X_{[3]}\tens  X_{[2]}\tens Y_{[1]} \ ,\cr
X_j\tens X[\underline e_i\times \underline e_j]\tens Y_i &=&
X_{[2]}\tens X_{[3]}\tens Y_{[1]}-X_{[1]}\tens X_{[3]}\tens Y_{[2]}-
X_{[3]}\tens X_{[2]}\tens Y_{[1]} \cr
&&+ X_{[1]}\tens X_{[2]}\tens Y_{[3]}+X_{[3]}\tens X_{[1]}\tens
Y_{[2]}-
X_{[2]}\tens X_{[1]}\tens Y_{[3]}\ ,
\end{eqnarray*}
and using this (\ref{vcbbvcdiblv}) simplifies to
\begin{eqnarray} \label{vcbbvcdiblvuu}
\phi
&=&  1\tens 1\tens 1+\hbar^2(Y[\underline e_i\times\underline
e_j]\tens  Y_j-
  X_j\tens X[\underline e_i\times \underline e_j] )\tens Y_i
/8+O(\hbar^3)\ .
\end{eqnarray}
Again expanding the summations, and assigning a name to part of
(\ref{vcbbvcdiblvuu}),
\begin{eqnarray*}
\psi&=&(Y[\underline e_i\times\underline e_j]\tens  Y_j-
  X_j\tens X[\underline e_i\times \underline e_j] )\tens Y_i \cr
  &=& (Y_{[3]}\tens Y_{[2]}-Y_{[2]}\tens Y_{[3]}
   - X_{[2]}\tens X_{[3]}  + X_{[3]}\tens X_{[2]})\tens Y_{[1]} \cr
   &&+  (Y_{[1]}\tens Y_{[3]}-Y_{[3]}\tens Y_{[1]}
   - X_{[3]}\tens X_{[1]}  + X_{[1]}\tens X_{[3]})\tens Y_{[2]} \cr
   && + (Y_{[2]}\tens Y_{[1]}-Y_{[1]}\tens Y_{[2]}
   - X_{[1]}\tens X_{[2]}  + X_{[2]}\tens X_{[1]})\tens Y_{[3]}
\end{eqnarray*}

\begin{lemma} If we take $\pi$ to be the reduction to $\tens_{C^\infty(S^2)}$, then
$\pi\psi=0$. 
\end{lemma}\proof This is a rather long calculation done with Mathematica. Details are omitted. \eproof 

This corresponds to the
multiplication of functions being associative to
$O(\hbar^2)$.  This critical fact justifies our choice of $G^{(2)}=(G^{(1)})^2/2$. Note that in this case we can use (\ref{forFinv}) to calculate
\begin{eqnarray*}
F &=& 1\tens 1 -\hbar\, G^{(1)}+\hbar^2\, (G^{(1)})^2/2+O(\hbar^3)\ 
\end{eqnarray*}
Note that this  would be consistent with $F^{-1}=e^{\hbar G^{(1)}}=e^{{\hbar\over 2}X_i\tens Y_i}$ but also note that $so(1,3)$ is not Abelian such an exponential form will not have $\phi=1\tens1\tens 1$, and indeed this is not true even at order $\hbar^2$. However, we see that when projected over $C^\infty(S^2)$ we do have that $\phi$ is effectively trivial at this order on functions; in other words the example demonstrates the hybrid set up of our paper at this order.

\subsection{The deformed algebra}
Following (\ref{prelim}), if $G^{(1)}=X_i\tens Y_i/2$ (summation suppressed)
and $G^{(2)}=(G^{(1)})^2/2$, then we have
\begin{eqnarray*}
f \bullet g &=& f\,g+ \hbar\, (X_i\la f)(Y_i\la g)/4+\hbar^2\, (X _iX_j\la
f)( Y_iY_j\la g)/8+O(\hbar^3)\ .
\end{eqnarray*}
It will be convenient to continue to use the coordiantes in
\ref{fycy}, in which case, using subscripts for partial
differentiation,
\begin{eqnarray*}
(X_i\la f)(Y_i\la g) &=&
((1-x^2)f_x-xyf_y)(-zg_y)+(-xyf_x+(1-y^2)f_y)zg_x \cr && +\,
(-xzf_x-yzf_y)(-yg_x+xg_y) \cr
&=& z\,(f_y\,g_x-f_x\,g_y)\ .
\end{eqnarray*}
With rather more work, we get the following formula:
\begin{eqnarray*}
(x^a\,y^b)\bullet(x^c\,y^d) &=& x^{a+c}\, y^{b+d} + \hbar\,z\,
x^{a+c-1}\, y^{b+d-1}(bc-ad)/2  \cr
&& + \hbar^2 \, x^{a+c-2}\, y^{b+d-2}\Big(
(bc(b-1)(c-1)+ad(a-1)(d-1)-2abcd) \cr
&&+y^2\,(a\,c - a^2\,c + b^2\,c - a\,c^2 - b^2\,c^2 + 2\,a\,b\,c\,d +
a\,d^2 -
   a^2\,d^2) \cr
   && +x^2(b\,c^2 - b^2\,c^2 + a^2\,d + b\,d - b^2\,d + 2\,a\,b\,c\,d
-
   a^2\,d^2 - b\,d^2) \cr
   &&+x^2\,y^2\left( a + b \right) \,\left( c + d \right) \,
   \left( 1 + a + b + c + d \right)
\Big)/8+O(\hbar^3)\ .
\end{eqnarray*}
The second order part of $f\bullet g$, evaluated at $x=y=0$, is one
eighth of
\begin{eqnarray*}
\frac{\partial^2 f}{\partial x^2}\frac{\partial^2 g}{\partial y^2}
+\frac{\partial^2 f}{\partial y^2}\frac{\partial^2 g}{\partial x^2}-
2\frac{\partial^2 f}{\partial x\partial y}\frac{\partial^2
g}{\partial x\partial y}-
\frac{\partial f}{\partial x}\frac{\partial g}{\partial x}
-\frac{\partial f}{\partial y}\frac{\partial g}{\partial y}\ .
\end{eqnarray*}
Unsing rotation invariance, and the fact that all the Christoffel
symbols
vanish at $x=y=0$, we see that this is
\begin{eqnarray*}
\omega^{ij}\,\omega^{kl}\,(\nabla_i f_{,k})\,(\nabla_j g_{,l})-
g^{ij}\,f_{,i}\,g_{,j}\ .
\end{eqnarray*}
This quantisation can be compared the Fedosov one for $S^2$ with the
above symplectic structure and symplectic connection.
We see that the second order part is not that given by the Fedosov
method,
as that does not have the $g^{ij}\,f_{,i}\,g_{,j}$ term.
Note that while Fedosov gives a prescription for a quantisation that
is associative on the functions to all orders from the symplectic
form and connection, this is not necessarily a unique quanitsation.
However we expect that our second order term $G^{(2)}$ may have to be
modified to allow extension to all orders in $\hbar$ as an
associative algebra, and it is not obvious that this could be done
within our existing 6 dimensional subalgebra of the
vector fields.

\section{Enveloping algebras $U(g)$ as cochain twists}  \label{jkascb}

As an important application of the ideas in Section~\ref{bvcihjvx},
we consider $M=\g^*$, the dual of a Lie algebra, equipped with its
standard Kirillov-Kostant Poisson structure $\{v,w\}=[v,w]$. Here
$S(\g)=\C[\g^*]$ i.e.\ we work with polynomial functions as generated
by $v\in\g$ viewed as linear functions on $\g^*$. 

\subsection{The cochain to lowest order.} 

To express these canonical data as induced by a cochain twist, we
seek suitable vector fields to define $G^{(1)}=\sum X\tens Y$. Some
natural vector fields are $\g$ itself acting by $\ad$ as mentioned
above, i.e. the vector fields for the coadjoint action on $\g^*$ from
a geometrical point of view. Then there is $\g^*$ acting by interior
product on $S(\g)$, which is to say usual differentiation on $\g^*$.
These classes of vector fields generated a sub-Lie algebra
$\CL=\g\rcross\g^*\subset\diff(\g^*)$ that turns out to be sufficient  to induce the
desired quantisation. Thus we take  $H=U({\frak g}\rcross{ \frak
g^*})=U({\frak g})\rcross S({\frak g^*})$ acting covariantly on
$A=S({\frak g})$.

Choose a basis $\{e_i\}$ in ${\frak g}$, and a dual basis $\{e^i\}$
in ${\frak g}^*$,
and set
\begin{eqnarray*}
F^{-1} &=& 1\tens 1 +\alpha \hbar\, e_i\tens e^i +\beta\hbar e^i\tens
e_i+O(\hbar^2)\ .
\end{eqnarray*}
Then
\[ \omega(\extd v,\extd w)=(\alpha-\beta)\sum e_i(v)
e^i(w)-e^i(v)e_i(w)=-2(\alpha-\beta)[v,w]\]
so we obtain the Kirillov-Kostant bracket with
$\beta-\alpha=\frac{1}{2}$. As a first consequence:

\begin{propos}
 $\g^*$ also has on it a canonical Poisson-compatible preconnection
\[ \varnabla_{v}\extd w=\frac{1}{2}\extd [v,w]\]
where $\hat v=\ad_v$ is the adjoint action viewed as a vector field
(in classical differential geometry, this is
a derivation on $S(\g)$). The curvature and torsion are 
\[ R(v,w)\extd z=-{1\over 4}[[v,w],z],\quad \<T(v,w),\extd z\>={1\over 2}[[v,w],z]\]
\end{propos} 
\proof This comes out of the construction  by using $F$ to deform the differential calculus, and the
leading order part $\sum X\tens Y$ found above. With hindsight one may check  independently using the axioms in \cite{BegMa:sem} that this is indeed a canonical Poisson-compatible
preconnection for the Kirillov-Kostant bracket. We  compute its curvature  
as
\[ R(v, w) (\extd z)
\,=\, ( \varnabla_{v} \varnabla_{w}
-\varnabla_{w} \varnabla_{v}-
\nabla_{\widehat{{\{v,w\}}}})(\extd z)=
\extd({[v,[w,z]]}-{[w,[v,z]]}-2{[[v,w],z]})/4\]
with the result stated in view of the Jacobi identity in the Lie algebra. Similarly 
\[ \<T( v,w),\extd z\>=\<\varnabla_{v}\hat w-\varnabla_{w}\hat v-[\hat v,\hat w],\extd z\>\\
=
\hat v(\hat w(z))-\<\hat w,\varnabla_{v}\extd z\>-(v\leftrightarrow w)-[[v,w],z] \]
with the result stated again on using the Jacobi identity. Incidentally, leaving out the $1/2$ gives a preconnection with zero curvature, but it is not Poisson-compatible. \eproof

Note that the most general translation-invariant  $\hat\nabla$ in this quantisation is given by the analysis of \cite{BegMa:sem} as of the form
\[ \hat\nabla_v\extd w=\frac{1}{2}\extd [v,w]+\extd \hat\Xi(v,w)\]
where $\hat\Xi:\cg\tens\cg\to\cg$ is some symmetric linear map. This follows from regarding $\cg^*$ as an Abelian Lie group and applying the theory in \cite[Sec. 4.1]{BegMa:sem}. The operations $L^*$ and $R^*$ translating differentials back to the origin are trivial so that $\hat\nabla_v\extd w=\extd\Xi(v,w)$ is defined by a map $\Xi$ with arbitrary symmetric part, which we denote $\hat\Xi$, and antisymmetric part the same as above.   On the other hand, if we further demand background `rotational' invariance in the sense of $\ad$-invariance under $\cg$ (which becomes covariance of the calculus under the quantum double $D(U(\cg))$ after quantisation) this corresponds to $\hat\Xi$ symmetric {\em and} ad-invariant. 

\begin{theorem} For all simple $\cg$ other than $sl_n$, $n>2$ the canonical preconnection in Proposition~5.1.1 is the only translation and $\cg$-invariant one on  $S(\cg)=\C[\cg^*]$. For $\cg=sl_n$, $n>2$ there is a 1-parameter moduli space of such covariant $\hat\nabla$ but they all have curvature. Hence for all simple $\cg$  any covariant differential calculus on $U_\hbar(\cg)$ with classical dimensions is necessarily nonassociative.  
\end{theorem}
\proof By the same arguments from invariant theory as in the proof of \cite[Theorem 4.20]{BegMa:sem} (but now in a different context), basically from Kostant's work,  there is no nonzero symmetric ad-invariant map $\hat\Xi:\cg\tens\cg\to \cg$ for $\cg$ simple other than for $sl_n$, $n>2$. Hence $\hat\Xi=0$ and $\hat\nabla$ has to be the one in Proposition~5.1.1. 

For $sl_n$, $n>2$ one has the possibility of a 1-parameter family via the invariant totally symmetric trilinear form viewed as the map $\hat\Xi$. In this case 
\[ R(v,w)\extd z=-{1\over 4}\extd [[v,w],z]+\extd\left(\hat\Xi(v,\hat\Xi(w,z))-\hat\Xi(w,\hat\Xi(v,z))\right)\]
since the terms linear in $\hat\Xi$ cancel using its $\ad$-invariance.   We have to show that there are always $v,w,z$ with the curvature expression nonzero.  To do this, note that the symmetric trilinear is a cubic polynomial on $sl_n$ which on $v\in sl_n$ has value $I(v,v,v)=I(v)$, say (e.g. for $sl_3$ we have $I(v)=\det(v)$). We can reconstruct the full trilinear from this by polarisation, e.g. 
\[ I(v,w,w)={1\over 6}(I(v+2w)-2I(v+w)+I(v))-I(w)\]
and we define $\hat\Xi(v,w)=I(v,w,e_i)e_j \kappa^{ij}$
where  $\kappa^{ij}$ is the inverse matrix of the Killing form (not necessarily normalised). We fix $v,w$ diagonal  (i.e. in the standard Cartan subalgebra of $sl_n$) and focus on
\[ R(v,w)\extd w=I(v,e_i,w)I(w,v,e_j)\kappa^{ij}-I(w,e_i,w)I(v,v,e_j)\kappa^{ij}.\]
We will show that this can be arranged to be non-zero. Note that if $t$ lies in the  Cartan subalgebra and $z\in sl_n$, then $\ad$-invariance  $I([t,v],w,z)+I(v,[t,w],z)+I(v,w,[t,z])=0$ means $I(v,w,[t,z])=0$. We conclude that if $z$ is a root vector then, $I(v,w,z)=0$ (since $[t,z]$ is a nonzero multiple of $z$). Hence $I(v,w,z)$ vanishes for all $z$ in the space spanned by the nonzero root vectors, which is to say the orthogonal complement of the Cartan with respect to the Killing form (it is the space of matrices in $sl_n$ with zero diagonal). Hence we let $\{e_a\}$ be a basis of the Cartan subalgebra completed to a basis of $sl_n$ taken from this complement. It means that we can compute $R(v,w)\extd w$ using only a sum over the $e_a,e_b$ in place of $e_i,e_j$ in the expression above. For $sl_3$ we take $t_1=e_{11}-e_{22}$, $t_2=e_{22}-e_{33}$ in the Cartan. Then $e_1=t_1$ and $e_2=t_1+{1\over 2}t_2={1\over 2}(e_{11}+e_{22})-e_{33}$ are a basis with $\kappa^{ab}={\rm diag}(1/2, 2/3)$. We also compute
\[ I_{111}\equiv I(t_1,t_1,t_1)=I_{222}\equiv I(t_2,t_2,t_2)=-1, \quad I_{112}\equiv I(t_1,t_1,t_2)={3\over 2},\quad I_{122}\equiv I(t_1,t_2,t_2)={5\over 6}\]
using the polarisation formula above. Hence setting $v=t_1$, $w=t_2$ we compute
\begin{eqnarray*}R(v,w)\extd w&=&{1\over 2}I(t_1,t_2,t_1)^2+{2\over 3}I(t_1,t_2,t_1+{1\over 2}t_2)^2\\
&&-{1\over 2}I(t_1,t_1,t_1)I(t_2,t_2,t_1)-{2\over 3}I(t_1,t_1,t_1+{1\over 2}t_2)I(t_2,t_2,t_1+{1\over 2}t_2)\\
&=&{1\over 2}I_{111}^2+{2\over 3}(I_{112}+{1\over 2}I_{122})^2-{1\over 2}I_{111}I_{122}-{2\over 3}(I_{111}+{1\over 2}I_{112})(I_{122}+{1\over 2}I_{222}) >0\end{eqnarray*}
for the values stated. This proves the result for $sl_3$. For $sl_n$ the trilinear is given by $I(v)=\sum_{i<j<k} v^iv^jv^k$ in terms of the diagonal entries of $v$ in the Cartan. We take the same $v=t_1,w=t_2$ as above but viewed in the standard way inside $sl_n$ rather than $sl_3$ and the $e_1,e_2$ completed to a diagonal basis for $\kappa$. Then the computation reduces to the same one as above for $sl_3$. \eproof

In summary, these results tells us that for simple $\cg$ we are going to necessarily have to work with nonassociative differentials, and for all Lie algebras $\cg$ the canonical `universal' choice (which is often the only choice) at the lowest order level is the one in Proposition~5.1.1. We therefore focus on this and have seen that it is indeed given by  a cochain twist at lowest order.  We next want to extend Proposition~5.1.1 to find $F,F^{-1}$ at least to order $O(\hbar^3)$. To do this we first 
look at the product of $U_\hbar(\g)$ on monomials. Computations have been done with MATHEMATICA.

\subsection{The Campbell-Baker-Hausdorff product}
We consider which $F$ induce not only the above semiclassical
data but the actual product of $U_\hbar(\g)$ as a star-product
quantisation of $S(g)$. Here $U_\hbar(\g)$ denotes the tensor algebra
on $\g$ with relations $vw-wv=\hbar [v,w]$ in terms of the Lie
bracket of $\g$. This is a deformation
of $S(\g)$ by the linear map  linear map $\varphi: S({\frak g}) \to
U_\hbar({\frak g})$ given by a sum over permutations
\begin{eqnarray*}
\varphi(v_1\dots v_n) &=& \frac{1}{n!}\sum_{\kappa\in S_n}
v_{\kappa(1)}\dots v_{\kappa(n)}\ .
\end{eqnarray*}
As explained in \cite{Gutt}, $\varphi$ is a 1-1 correspondence, and
we define a deformed multiplication $\bullet$ on $S({\frak g})$ by
$\underline v\bullet\underline w=\varphi^{-1}(\varphi(\underline v).
\varphi(\underline w))$. As examples,
\begin{eqnarray*}
\varphi(v).\varphi(w) &=& (vw+wv)/2 + \hbar \,[v,w]/2 \cr
&=& \varphi(vw) + \hbar\,\varphi([v,w])/2\ , \cr
\varphi(v^2).\varphi(w) &=& (v^2w+vwv+wv^2)/3+\hbar\,(
2\,v[v,w]+[v,w]v)/3 \cr
&=&
\varphi(v^2w)+\hbar\,\varphi(v[v,w])/2+\hbar^2\,\varphi([v[v,w]])/6\ .
\end{eqnarray*}
This is related to the CBH formula  as follows: since $e^v$ in $S(\g)$
maps under $\phi$ to $e^v$ in $U_\hbar(\g)$ (similarly for any power
series in $v$) we have
\[ e^v\bullet
e^w=\phi^{-1}(\phi(e^v)\phi(e^w))=\phi^{-1}(e^ve^w)=\phi^{-1}(e^{C_\hbar(v,w)})=e^{C_\hbar(v,w)}\]
where $C_\hbar(v,w)$ is the CBH power series for the product of two
exponentials in $U(\g)$ with the insertion of powers of $\hbar$ for
each commutator in the Lie algebra.

\begin{lemma} \label{vchcsavb}
For $v_1\dots v_n$ a symmetric product of elements of ${\mathfrak
g}$, and
$w\in {\mathfrak g}$,
\begin{eqnarray*}
w\bullet v_1\dots v_n &=& \varphi(wv_1\dots v_n)-\hbar\,n\,
\varphi([v_1,w]v_2\dots v_n)/2 \cr
&& +\, \hbar^2\,n(n-1)\varphi([v_1,[v_2,w]]v_3\dots
v_n)/12+O(\hbar^3)\ .
\end{eqnarray*}
\end{lemma}
\noindent {\bf Proof:}\quad  Set
\begin{eqnarray*}
0\le i\le n \qquad a^0_i &=& v_1\dots v_i w\dots v_n\ , \cr
0\le i\le n-1 \qquad a^1_i &=& v_1\dots v_i [v_{i+1},w]\dots v_n\ ,\cr
0\le i\le n-2 \qquad a^2_i &=& v_1\dots v_i
[v_{i+1},[v_{i+2},w]]\dots v_n\ .
\end{eqnarray*}
Then in $U_\hbar({\mathfrak g})$ we have
\begin{eqnarray*}
a^m_i-a^m_{i+1} &=& -\hbar\, a^{m+1}_{i}\ ,
\end{eqnarray*}
and this gives, for $i>0$,
\begin{eqnarray} \label{vbcus}
a^m_i &=& a^m_0 + \hbar(a^{m+1}_0+\dots+a^{m+1}_{i-1})\ .
\end{eqnarray}
 From this we get
\begin{eqnarray} \label{vbcashu}
a^0_0+\dots+a^0_n &=& (n+1)a^0_0 + \hbar
(n\,a^{1}_0+(n-1)a^1_1+\dots+a^{1}_{n-1}) \cr
&=& (n+1)a^0_0  + \hbar(n+1)(a^1_0+\dots+a^1_{n-1})/2 \cr
  && + \hbar
((n-1)\,a^{1}_0+(n-3)a^1_1+\dots+(1-n)a^{1}_{n-1})/2\ .
\end{eqnarray}
Now we use (\ref{vbcus}) again to get
\begin{eqnarray*}
\sum_{j:0\le j\le n-1}(n-1-2j)\,a^1_j &=& \hbar \sum_{j:0\le j\le
n-1}(n-1-2j)\,
\sum_{i:0\le i\le j-1} a^2_i \cr
&=& \hbar \sum_{i:0\le i\le n-2}a^2_i
\sum_{j:i+1\le j\le n-1}(n-1-2j)\cr
&=& -\,\hbar \sum_{i:0\le i\le n-2}a^2_i
(i+1)(n-1-i)\ ,
\end{eqnarray*}
and from this (\ref{vbcashu}) becomes the following, where we use
(\ref{vbcus}) again to get the last equality:
\begin{eqnarray}
a^0_0+\dots+a^0_n
&=& (n+1)a^0_0  + \hbar(n+1)(a^1_0+\dots+a^1_{n-1})/2 \cr
  && -\, \frac{ \hbar^2}{2}
  \sum_{i:0\le i\le n-2}a^2_i   (i+1)(n-1-i)   \cr
  &=& (n+1)a^0_0  + \hbar(n+1)(a^1_0+\dots+a^1_{n-1})/2 \cr
  && -\, \frac{n(n+1) \hbar^2}{12}(a^2_0+\dots+a^2_{n-2})+O(\hbar^3)
\end{eqnarray} 
as required. \eproof

\begin{lemma}
For $v_1\dots v_n$ and
$w_0\dots w_m$ symmetric products of elements of ${\mathfrak g}$,
\begin{eqnarray*}
w_0\bullet  \varphi(w_1\dots w_m v_1\dots v_n) &=&
\varphi(w_0w_1\dots w_m v_1\dots v_n)\cr
  && -\hbar\,n\,\varphi([v_1,w_0]w_1\dots w_m v_2\dots v_n)/2\cr
  &&+\hbar^2\,\varphi\big(  n(n-1) [v_1,[v_2,w_0]]w_1\dots w_m
v_3\dots v_n\cr
&& + nm[w_1,[v_1,w_0]]w_2\dots w_m v_2\dots v_n\big)/12+O(\hbar^3) \cr
w_0\bullet  \varphi([v_1,w_1]w_2\dots w_m v_2\dots v_n) &=&
\varphi([v_1,w_0]w_1\dots w_m v_2\dots v_n) \cr
&& -\hbar\, \varphi\big( [w_0,[w_1,v_1]]w_2\dots w_m v_2\dots v_n \cr
&&+(n-1)[v_1,w_0][v_2,w_1]w_2\dots w_m v_3\dots
v_n\big)/2+O(\hbar^2)\ .
\end{eqnarray*}
\end{lemma}
\noindent {\bf Proof:}\quad Using \ref{vchcsavb}
and being careful about counting permutations, we get, for $u\in
{\mathfrak g}$,
to $O(\hbar^3)$,
\begin{eqnarray*}
u\bullet  \varphi(w_1\dots w_m v_1\dots v_n) &=&
\varphi(uw_1\dots w_m v_1\dots v_n)\cr
  && -\hbar\,(n+m)\,n\,\varphi([v_1,u]w_1\dots w_m v_2\dots
v_n)/(2(n+m))\cr
  && -\hbar\,(n+m)\,m\,\varphi([w_1,u]w_2\dots w_m v_1\dots
v_n)/(2(n+m)) \cr
  &&+\hbar^2\,(n+m)(n+m-1)(1/12)\,\varphi\big( \cr
&& n(n-1) [v_1,[v_2,u]]w_1\dots w_m v_3\dots v_n/((n+m)(n+m-1))\cr
&& + nm [v_1,[w_1,u]]w_2\dots w_m v_2\dots v_n/((n+m)(n+m-1))\cr
&& + nm[w_1,[v_1,u]]w_2\dots w_m v_2\dots v_n/((n+m)(n+m-1))\cr
&& + m(m-1)[w_1,[w_2,u]]w_3\dots w_m v_1\dots
v_n/((n+m)(n+m-1))\big)\ .
\end{eqnarray*}
Putting $u=w_0$ and supposing that $w_0\dots w_m$ is symmetrised,
this reduces to
the first part of the statement.
Next, to $O(\hbar^2)$,
\begin{eqnarray*}
u\bullet  \varphi([v_1,w_1]w_2\dots w_m v_2\dots v_n) &=&
\varphi(u[v_1,w_1]w_2\dots w_m v_2\dots v_n) \cr
&& -\hbar(n+m-1)(1/2)\varphi\big( \cr
&&[[v_1,w_1],u]w_2\dots w_m v_2\dots v_n/(n+m-1) \cr
&&+(m-1)[w_2,u][v_1,w_1]w_3\dots w_m v_2\dots v_n/(n+m-1) \cr
&&+(n-1)[v_2,u][v_1,w_1]w_2\dots w_m v_3\dots v_n/(n+m-1) \big)\ .
\end{eqnarray*}
If $w_0\dots w_m$ is symmetrised, this reduces to the second part of
the statement. \eproof

\begin{propos} \label{jhgsadfv}
For symmetric $w_1\dots w_m$
and $v_1\dots v_n$:
\begin{eqnarray*}
w_1\dots w_m\bullet v_1\dots v_n &=&  \varphi(w_1\dots w_m v_1\dots
v_n)\cr &&  -
\hbar\, mn\,\varphi([v_1,w_1]w_2\dots w_m v_2\dots v_n)/2 \cr
&&+\hbar^2\,n(n-1)m(m-1)\,\varphi([v_1,w_1][v_2,w_2]w_3\dots w_m
v_3\dots v_n)/8  \cr
&&+\hbar^2\,n(n-1)m\,\varphi([v_1,[v_2,w_1]]w_2\dots w_m v_3\dots
v_n)/12  \cr
&&+\hbar^2\,m(m-1)n\,\varphi([w_1,[w_2,v_1]]w_3\dots w_m v_2\dots
v_n) /12
+O(\hbar^3)\ .
\end{eqnarray*}
\end{propos}
\noindent {\bf Proof:}\quad By induction on $m$.
We will suppose that, for fixed $n$ and symmetric $w_1\dots w_m$,
\begin{eqnarray*}
w_1\dots w_m\bullet v_1\dots v_n &=&  \varphi(w_1\dots w_m v_1\dots
v_n) -
\hbar\, \alpha_m\,\varphi([v_1,w_1]w_2\dots w_m v_2\dots v_n) \cr
&&+\hbar^2\big(\beta_m\,\varphi([v_1,w_1][v_2,w_2]w_3\dots w_m
v_3\dots v_n)  \cr
&&+\gamma_m\,\varphi([v_1,[v_2,w_1]]w_2\dots w_m v_3\dots v_n)  \cr
&&+\delta_m\,\varphi([w_1,[w_2,v_1]]w_3\dots w_m v_2\dots v_n)
\big)+O(\hbar^3)\ .
\end{eqnarray*}
Using associativity of the $\bullet$ product,
\begin{eqnarray*}
w_0\dots w_m\bullet v_1\dots v_n &=&
\varphi(w_0w_1\dots w_m v_1\dots v_n)\cr
  && -\hbar\,n\,\varphi([v_1,w_0]w_1\dots w_m v_2\dots v_n)/2\cr
  &&+\hbar^2\,\varphi\big(  n(n-1) [v_1,[v_2,w_0]]w_1\dots w_m
v_3\dots v_n\cr
&& - nm[w_0,[w_1,v_1]]w_2\dots w_m v_2\dots v_n\big)/12\cr
  &&
-\hbar\,\alpha_m\, \varphi([v_1,w_0]w_1\dots w_m v_2\dots v_n) \cr
&& +\hbar^2\,\alpha_m\, \varphi\big( [w_0,[w_1,v_1]]w_2\dots w_m
v_2\dots v_n \cr
&&+(n-1)[v_1,w_0][v_2,w_1]w_2\dots w_m v_3\dots v_n\big)/2 \cr
&&+\hbar^2\big(\beta_m\,\varphi([v_1,w_0][v_2,w_1]w_2\dots w_m
v_3\dots v_n)  \cr
&&+\gamma_m\,\varphi([v_1,[v_2,w_0]]w_1\dots w_m v_3\dots v_n)  \cr
&&+\delta_m\,\varphi([w_0,[w_1,v_1]]w_2\dots w_m v_2\dots v_n)
\big)+O(\hbar^3)\ .
\end{eqnarray*}
This gives the recursive equations and initial conditions
\begin{eqnarray*}
\alpha_{m+1} &=& \alpha_m +n/2\ ,\quad \alpha_1\,=\,n/2\ ,\cr
\beta_{m+1} &=& \beta_m +\alpha_m\,(n-1)/2\ ,\quad \beta_1\,=\,0\ ,\cr
\gamma_{m+1} &=& \gamma_m + n(n-1)/12\ ,\quad \gamma_1\,=\,n(n-1)/12\
,\cr
\delta_{m+1} &=& \delta_m -nm/12+\alpha_m/2\ ,\quad\delta_1\,=\,0 \ .
\end{eqnarray*}
{}From this we get $\alpha_m=nm/2$, $\gamma_m=n(n-1)m/12$,
$\delta_m=m(m-1)n/12$ and $\beta_m=n(n-1)m(m-1)/8$. \eproof

\subsection{Cochain for the deformed product of $S({\frak g})$}
Choose a dual basis $e_i\in {\frak g}$ and $e^i\in {\frak g}^*$. Let
${\frak g}$
act on $S({\frak g})$ by the adjoint, and ${\frak g}^*$ act by
evaluation.
Set $Q_1=e_i\tens e^i$ and $Q_2=e^i\tens e_i$, and let
$\mu$ stand for multiplication. Then for symmetric $\underline w =
w_1\dots w_m$ and $\underline v=v_1\dots v_n$:
\begin{eqnarray*}
\mu(e_i\tens e^i)(\underline w\tens \underline v) &=&
  n\,[v_1,w_1\dots w_m] v_2\dots v_n \cr
  &=&  nm\,[v_1,w_1]w_2\dots w_m v_2\dots v_n \cr
\mu(e^i\tens e_i)(\underline w\tens \underline v) &=&
  m\,w_2\dots w_m [w_1,v_1\dots v_n] \cr
   &=&  -nm\,[v_1,w_1]w_2\dots w_m v_2\dots v_n  \ .
\end{eqnarray*}
We have quadratic terms of the form $Q_1^2$, $Q_2^2$ and $:\!Q_1Q_2\!:$
(with $::$ denoting a normal ordering
with elements of ${\frak g}^*$ being put on the right), which are
respectively
\begin{eqnarray*}
\mu(e_ie_j\tens e^ie^j)(\underline w\tens \underline v) &=&
n(n-1)\, [v_1,[v_2,w_1\dots w_m]] v_3\dots v_n\cr
&=&
n(n-1)m\, [v_2,[v_1,w_1]w_2\dots w_m] v_3\dots v_n\cr
&=&  n(n-1)m\, [v_2,[v_1,w_1]]w_2\dots w_m v_3\dots v_n \cr &&+\,
n(n-1)m(m-1)\, [v_1,w_1][v_2,w_2]w_3\dots w_m v_3\dots v_n\ ,  \cr
\mu(e^ie^j\tens e_ie_j)(\underline w\tens \underline v) &=&
m(m-1)\, w_3\dots w_m [w_1,[w_2,v_1\dots v_n]]  \cr
&=& nm(m-1)\, w_3\dots w_m [w_1,[w_2,v_1]v_2\dots v_n]  \cr
&=& nm(m-1)\, w_3\dots w_m [w_1,[w_2,v_1]]v_2\dots v_n  \cr
&& +\,n(n-1)m(m-1)\, w_3\dots w_m [w_2,v_2][w_1,v_1]v_3\dots v_n \
,\cr
\mu(e_ie^j\tens e_je^i)(\underline w\tens \underline v) &=& n [v_1,
e^j(\underline w)]
. e_j(v_2\dots v_n) \cr
&=& nm\,[v_1, w_2\dots w_m]. [w_1,v_2\dots v_n] \cr
&=& -n(n-1)m(m-1)\, [v_1,w_1][v_2,w_2]w_3\dots w_m v_3\dots v_n\ .
\end{eqnarray*}
If we set
\begin{eqnarray} \label{kjhsdvc}
F^{-1} &=& 1\tens 1 +\hbar(\alpha \, Q_1+\beta\,Q_2)+\hbar^2(2\,Q_1^2+
2\, Q_2^2+\! :\!Q_1 Q_2\!:)/24+O(\hbar^3)\ ,
\end{eqnarray}
where $\alpha-\beta \,=\, -1/2$, then we recover the CBH product to
$O(\hbar^3)$. Note that we could add any multiple of $Q_1Q_2+Q_2^2$ or 
$Q_2Q_1+Q_1^2$ to $G^{(2)}$ and still get the same product to
$O(\hbar^3)$. Also we have the equation
\begin{eqnarray*}
\mu(:\!Q_1Q_2\!:^R+:\!Q_1Q_2\!:+Q_1^2+Q_2^2)(\underline w\tens \underline v) &=&
nm\,w_2\dots w_m v_2\dots v_n.e^i([v_1,[w_1,e_i]])\ ,
\end{eqnarray*}
where $:\!Q_1Q_2\!:^R=e^je_i\tens e^ie_j$ is the reversed normal order. Now
$e^i([v,[w,e_i]])$ (summed over $i$) is the trace of ${\rm ad}_v{\rm ad}_w$, that is
$-\<v,w\>$, where $\<,\>$ is the killing form. If we set $\<e_i,e_j\>=\kappa_{ij}$, then
\begin{eqnarray*}
\mu(:\!Q_1Q_2\!:^R+:\!Q_1Q_2\!:+\,Q_1^2+Q_2^2+\kappa_{ij}\,e^i\tens e^j
)(\underline w\tens \underline v) &=&0\ .
\end{eqnarray*}

\subsection{Improved cochain for the deformed coproduct on $S({\frak g}^*)$}

Here we consider a more general covariant ansatz but show that a further requirement relating to the coproduct of $S({\cg})$ again leads to a unique answer.  Thus, we can write a more general expression for $G^{(2)}$ as
\begin{eqnarray} \label{g2new}
G^{(2)} &=& (2\,Q_1^2+
2\, Q_2^2+\! :\!Q_1 Q_2\!:)/24 + \gamma (Q_1+Q_2)Q_1+
\delta (Q_1+Q_2)Q_2 \cr
&&\,+ \, \zeta(:\!Q_1Q_2\!:^R+:\!Q_1Q_2\!:+\,Q_1^2+Q_2^2+\kappa_{ij}\,e^i\tens e^j)\ ,
\end{eqnarray}
where $\gamma,\delta$ and $\zeta$ are constants,
and we shall use this instead of the $\hbar^2$ term in (\ref{kjhsdvc}). By the 
above discussion, this still gives the correct deformed product on $S({\frak g})$. 
We use $G^{(1)}=\alpha\,Q_1+\beta\,Q_2$. Then from (\ref{forFinv}) the deformed
coproduct is given by
\begin{eqnarray} \label{ksujkcv}
\Delta_F(x)\,=\, \Delta(x)+\hbar\,[\Delta(x),G^{(1)}]+\hbar^2\,([\Delta(x),G^{(2)}]-
G^{(1)}[\Delta(x),G^{(1)}])+O(\hbar^3)\ .
\end{eqnarray}
For $x\in {\frak g}^*$ we have $\Delta(x)=x\tens 1+1\tens x$, and using the
fact that elements of  ${\frak g}^*$ commute, we find
\begin{eqnarray*}
[\Delta(x),G^{(1)}] \,=\, \alpha\,[x,e_i]\tens e^i+\beta\, e^i\tens [x,e_i]\ .
\end{eqnarray*}
The coefficient of $\hbar^2$ in (\ref{ksujkcv}) is
\begin{eqnarray*}
&&\Big( \gamma[x,e_i]e_j+(\gamma-\alpha^2)e_j[x,e_i]+(\frac1{12}+\zeta)
[[x,e_i],e_j]\Big)\tens e^ie^j \cr
&&+\, e^ie^j\tens \Big(
\delta[x,e_i]e_j+(\delta-\beta^2)e_j[x,e_i]+(\frac1{12}+\zeta) [[x,e_i],e_j]
\Big) \cr
&&+\,e^i[x,e_j]\tens\Big(
(\gamma+\frac1{24}-\alpha\beta+\zeta)e_ie^j+(\delta+\zeta) e^je_i
\Big) \cr
&&+\,\Big(
(\gamma+\zeta) e^i e_j+(\delta+\frac1{24}-\alpha\beta+\zeta)e_je^i
\Big) \tens [x,e_i]e^j\ .
\end{eqnarray*}
To ensure that this is in $S({\frak g}^*)\tens S({\frak g}^*)$ we require that
$\gamma=\alpha^2/2$, $\delta=\beta^2/2$ and $(\alpha-\beta)^2=-4\zeta-1/12$. 
We already have $\alpha-\beta=-1/2$, so we get $\zeta=-1/12$. 
Now the 
coefficient of $\hbar^2$ in (\ref{ksujkcv}) is:
\begin{eqnarray*}\Delta_F(x)&=&x\tens 1+1\tens x+\hbar(\alpha[x,e_i]\tens e^i+\beta e^i\tens[x,e_i]\\
&&+\hbar^2 \alpha^2\,[[x,e_i],e_j] \tens e^ie^j/2 + e^ie^j\tens \beta^2 [[x,e_i],e_j]/2 \\
&&+\,(\beta^2/2-1/12) \,e^i[x,e_j]\tens[e^j,e_i] +(\alpha^2/2-1/12)\, [e^i, e_j] \tens [x,e_i]e^j\ .
\end{eqnarray*}

Also putting these values into the general form of $F^{-1}$ we obtain:

 \begin{theorem} The cochain 
\begin{eqnarray*} 
F^{-1}&=&1\tens 1+\hbar(\alpha Q_1+\beta Q_2)\\
&&+\hbar^2({\alpha^2\over 2}(Q_1+Q_2)Q_1+{\beta^2\over 2}(Q_1+Q_2)Q_2-{1\over 24}(:Q_1Q_2:+2:Q_1Q_2:^R+2\kappa_{ij}e^i\tens e^j))+O(\hbar^3)
\end{eqnarray*}
for $\alpha-\beta=-{1\over 2}$ reproduces the product of $U_\hbar(\cg)$ to the relevant order and 
has the property that $\Delta_F(S(\cg^*))\subset S(\cg^*)\tens S(\cg^*)$ to the relevant order.  
\end{theorem}
 
It appears that the up to the choice of how $-{1\over 2}$ is split between $\alpha$ and $\beta$, the cochain $F$ is determined at higher orders by the properties in the theorem to hold for any Lie algebra and the requirement of having a $\cg$-invariant form.  Although we will not give a formal proof, let us explain the underlying reason here. First we note that
locally near the identity we may identify the Lie algebra with the connected and simply connected Lie group $G$ associated to it, i.e. $S(\cg^*)_F\subseteq \C_{loc}[G]$, where the latter denotes functions defined near the identity. This is by
\[ \Theta:S(\cg^*)\supset \cg^*\owns x\mapsto f_x(e^{\hbar v})=\hbar^{-1}\<x,v\>\]
(i.e. the generators appear as logarithmic coordinates on the Lie group). Let us now show that
\[ ((\Theta\tens\Theta)\Delta_F x)(e^v,e^w)=(\Delta_{\C_{loc}[G]}\Theta(x))(e^v,e^w),\]
i.e. the identification is indeed as Hopf algebras to the relevant order. The right hand side here is $\Theta(x)(e^ve^w)=\Theta(x)(e^{C(v,w)})=\hbar^{-1}\<x,C(v,w)\>$ where $C(v,w)$ is the Campbell-Baker-Hausdorf series. Meanwhile, evaluating the coproduct $\Delta_F$ we have
\begin{eqnarray*}
&&\kern-20pt ((\Theta\tens\Theta)\Delta_F x)(e^v,e^w)=\hbar^{-1}(\<x,v\>+\<x,w\>)+\hbar^{-1}\alpha\<[x,e_i],v\>\<e^i,w\>+\hbar^{-1}\beta\<e^i,v\>\<[x,e_i],w\>\\
&&+\hbar^{-1}{\alpha^2\over 2}\<v,[[x,e_i],e_j]\>\<e^i,w\>\<e^j,w\>+\hbar^{-1}{\beta^2\over 2}\<w,e^i\>\<w,e^j\>\<v,[[x,e_i],e_j]\>\\
&&+\hbar^{-1}({\alpha^2\over 2}-{1\over 12})\<v,e^i\>\<v,[x,e_j]\>\<w,[e^j,e_i]\>
+\hbar^{-1}({\beta^2\over 2}-{1\over 12})\<v,e^i\><v,e^j\>\<w,[[x,e_i],e_j]\>+\cdots\\
&&=\hbar^{-1}\<x,v+w+[w,v](\alpha-\beta)+{\alpha^2\over 2}[w,[w,v]]+{\beta^2\over 2}[v,[v,w]]+({\alpha^2\over 2}-{1\over 12})[[w,v],w]\\
&&\quad+({\beta^2\over 2}-{1\over 12})[[v,w],v]\>+\cdots \\
&&=\hbar^{-1}\<x,v+w+{1\over 2}[v,w]+{1\over 12}([v,[v,w]]+[[v,w],w])+\cdots \>\end{eqnarray*}
for all $x\in\cg^*$. We used $\<[x,v],w\>=\<x\ra v,w\>=\<x,v\la w\>=\<x,[v,w]\>$ and similarly for repeated commutators. Note that the increasing powers of $\hbar^{-1}$ with each evalution exactly match the increasing powers of $\hbar$ in the powerseries coming from $F$. Hence the twisted coproduct reproduces  the
Campbell-Baker-Hausdorf series to low degree in its expansion. It is clear that this requirement and that we continue to reproduce the product of $U_\hbar(\cg)$ can be used to determine a universal formula for $F,F^{-1}$, though again it would be beyond our scope to provide this here.

\subsection{The Duflo map}

The Duflo map provides an independent check of our formula in Theorem~5.4.1 and gives some idea of  
the structure of $F^{-1}$ at all orders. We recall \cite{Du} that there is an invertible operator $D$ on $S(\cg)$ defined by
\[ D=e^{\sum_{k=1}^\infty \alpha_{2k}\del_{\Tr_{2k}}};\quad \sum_{k=1}^\infty \alpha_{2k}t^{2k}\equiv{1\over 2}\ln\left({\sinh(t/2)\over t/2}\right)\]
where $\del_{\Tr_{2k}}$ is the differential operator on $S(\cg)$ given by the action of the element $\Tr_{2k}=\Tr_\cg((\ad_\cdot)^{2k})\in S(\cg^*)$. The lowest order part is
\[ \alpha_{2}={1\over 48},\quad \alpha_{4}=-{1\over 5760},\quad \del_{\Tr_2}=-\kappa_{ij}e^ie^j\]
in our conventions above.  Duflo's theorem is that when restricted to the ad-invariant subalgebra $S(\cg)^\cg$ the map $\varphi\circ D$ is an isomorphism of this with the centre $Z(U(\cg))$.

Now let $F_{\rm red}^{-1}\in U(\cg^*)^{\tens 2}$ be the effective $F^{-1}$ when acting on invariants $S(\cg)^\cg\tens S(\cg)^\cg$. This is given by normal ordering $F^{-1}$ so that all terms have all elements of $\cg$ to the right of all elements of $\cg^*$ (what we called $:\ :^R$ above). Then project $\cg$ to zero in the result because by definition it act by zero on invariant elements. The result is some power-series in $U(\cg^*)^{\tens 2}=S(\cg^*)^{\tens 2}$ since $\cg^*$ is being regarded as an Abelian Lie algebra. 

\begin{propos} At least to $O(\hbar^3)$, $F^{-1}_{\rm red}$ for the cochain in Theorem~5.4.1 is a coboundary in the sense of \cite[Chapter 2.3]{Ma:book} of the Duflo element,  i.e.,
\[ F^{-1}_{\rm red}=(\Delta \gamma )(\gamma^{-1}\tens \gamma^{-1}),\quad \gamma=e^{\sum_{k=1}^\infty \alpha_{2k}\hbar^{2k}{\Tr_{2k}}}\]
\end{propos}

Here $\gamma$ viewed as an operator acting on $S(\cg)$ is just $D$ in Duflo's theorem after explicitly introducing the deformation scaling parameter. The coproduct $\Delta$ is that of $U(\cg^*)$. We expect this result to hold to all orders because $F^{-1}_{\rm red}$ a coboundary of some cochain $\gamma$ implies that
\[f\bullet g=\mu(F^{-1}_{\rm red}.(f\tens g))=\gamma.\mu(\gamma^{-1}.f\tens\gamma^{-1}.g))=D(\mu(D^{-1}f\tens D^{-1}g))\]
for all $f,g\in S(\cg)^\cg$, where $D$ denotes $\gamma$ acting on $S(\cg)$. We used in the first equality that $U(\cg^*)$ acts covariant on $S(\cg)$ with its initial product $\mu$ and hence we can move the action of $\Delta\gamma$ to the left as the action of $\gamma$. This means that the modified product restricted to invariant elements is an isomorphism of algebras (this is the meaning of $F^{-1}_{\rm red}$ being a coboundary as explained in \cite{Ma:book}). In the light of Duflo's theorem we explect $F^{-1}_{\rm red}$ therefore to be a coboundary of an invertible element $\gamma$ whose action is the same as the operator $D$ in Duflo's theorem. This leads to the statement of the proposition. 

We now verify the proposition to the order $O(\hbar^3)$ available to us. In the expression in Theorem~5.4.1
all the terms have some $e_i$ already to the right and therefore fail to contribute, except $:Q_1Q_2:$ and the $\kappa_{ij}e^i\tens e^j$ terms. We write the former using
\[ e_je^i=e^ie_j+(e^i\ra e_j)=e^ie_j+f_{jki}e^k\]
where $[e_i,e_j]=f_{ijk}e_k$ defines the structure constants. Then
\[ :Q_1Q_2:=e_ie^j\tens e_je^i\sim f_{jki}f_{imj}e^k\tens e^m\sim f_{kji}f_{mij}e^k\tens e^m\sim -\kappa_{ij}e^i\tens e^j\]
discarding terms acting trivially on invariant elements. As a result we have
\[ F^{-1}_{\rm red}=1\tens 1- {\hbar^2\over 24}\kappa_{ij}e^i\tens e^j+O(\hbar ^3)\]
(in fact the next term should be $O(\hbar^4)$). On the other hand from the above
\[ \gamma=e^{-{\hbar^2\over 48}c+O(\hbar^4)};\quad c=\kappa_{ij}e^ie^j,\quad (\Delta\gamma)(\gamma^{-1}\tens\gamma^{-1})=e^{-{\hbar^2\over 48}(\Delta c-c\tens 1-1\tens c)+O(\hbar^4)}=F^{-1}_{\rm red}\]
 to lowest order. The same proposition   would provide a check to  all orders of any cochain found. At the moment we have provided a check of our order $O(\hbar^3)$ result.

\subsection{Example:
noncommutative Minkowski space as cochain twist.} Here we verify (\ref{kjhsdvc})  for the algebra
 $[t,x_i]=\hbar x_i$ which has been proposed as noncommutative spacetime (the so-called bicrossproduct model). For convenience we take only one $x=x_i$ rather than $i=1,2,3$ for spacetime, however the structure is exactly similar. In this case we  exhibit a candidate for $F^{-1}$ to the next order, i.e.  up to $O(\hbar^4)$.  
  
 In this model there is a representation of the algebra in terms of $2\times 2$
matrices, as
  \begin{eqnarray*}
t \mapsto \left(\begin{array}{cc}\hbar/2 & 0 \\0 &
-\hbar/2\end{array}\right)\ ,
\quad x\mapsto \left(\begin{array}{cc}0 & 1 \\0 &
0\end{array}\right)\ .
\end{eqnarray*}
These matrices can be exponentiated to give
\begin{eqnarray*}
\exp(pt+qx) &=& \left(\begin{array}{cc}e^{\frac{\hbar\,p}{2}} &
\frac{\left( -1 +
        e^{\hbar\,p} \right) \,q}{e^{\frac{\hbar\,p}{2}}\,\hbar\,p}
\\0 & e^{\frac{-\left( \hbar\,p \right) }{2}}\end{array}\right)\ .
\end{eqnarray*}
A little matrix multiplication shows that
\begin{eqnarray*}
\exp(pt+qx).\exp(rt+sx) &=& \exp\big((p+r)t+\frac{\left( p + r
\right) \,
     \left( \left( -1 + e^{\hbar\,p} \right) \,q\,r +
       e^{\hbar\,p}\,\left( -1 + e^{\hbar\,r} \right) \,p\,s
       \right) }{\left( -1 +
       e^{\hbar\,\left( p + r \right) } \right) \,p\,r}\, x\big)\ ,
\end{eqnarray*}
so in this case we have a closed form for the CBH formula.

Further calculation with this algebra gives
\begin{eqnarray*}
\varphi(x^nt^m) &=& x^n\Big(t^m+\frac{nm}{2}\hbar\, t^{m-1}+
\frac{n(3n+1)m(m-1)}{24}\hbar^2 \,t^{m-2} \cr
&&+\frac{n^2(n+1)m(m-1)(m-2)}{48}\hbar^3\, t^{m-3}+O(\hbar^4)\Big)\ .
\end{eqnarray*}
 From this we can calculate
\begin{eqnarray*}
\varphi^{-1}(x^nt^m) &=& x^n\Big(t^m-\frac{nm}{2}\hbar\, t^{m-1}
+\frac{n(3n-1)m(m-1)}{24}\hbar^2 \,t^{m-2} \cr
&&+\frac{n^2(1-n)m(m-1)(m-2)}{48}\hbar^3\, t^{m-3}+O(\hbar^4)\Big)\ .
\end{eqnarray*}
If we combine this with the following formula for multiplication in
$U_\hbar$,
\begin{eqnarray*}
t^s\,x^r\,=\,\sum_{p=0}^s C^s_p\,(r\,\hbar)^p\, x^r\,t^{s-p}\ ,
\end{eqnarray*}
(where $C^s_p$ is a binomial coefficient) we get the formula
\begin{eqnarray} \label{ujzvsc}
(x^n\,t^m)\bullet(x^r\,t^s) &=&
x^{n+r}\,t^{m+s}+\frac{\hbar}2(mr-ns)\,x^{n+r}\,t^{m+s-1}
\cr && +\, \frac{\hbar^2}{24}\Big(mr - m^2r - 3mr^2 + 3m^2r^2 + ns -
   2mns - 3n^2s - 2mrs - 6mnrs  \cr
&& -
   ns^2 + 3n^2s^2\Big)\,x^{n+r}\,t^{m+s-2} + \frac{\hbar^3}{48} \Big(
-2mr^2 + 3m^2r^2 - m^3r^2 + 2mr^3  \cr && -
   3m^2r^3 + m^3r^3 + 2n^2s - 2mn^2s -
   2n^3s - m^2nrs - 3mn^2rs +
   2mr^2s  \cr &&  - 2m^2r^2s + 3mnr^2s -
   3m^2nr^2s - 3n^2s^2 + 2mn^2s^2 +
   3n^3s^2 + mnrs^2 \cr && + 3mn^2rs^2 +
   n^2s^3 - n^3s^3\Big)\,x^{n+r}\,t^{m+s-3} + O(\hbar^4)\ .
\end{eqnarray}
On $S({\frak g})$, $\ad_t$ is identified with $x\,\frac{\extd}{\extd
x}$,
and $\ad_x$ is identified with $-x\,\frac{\extd}{\extd t}$. We take
the dual
basis $\hat t,\hat x\in {\frak g}^*$. Then on $S({\frak g})$, $\hat
t$ is identified
with $\frac{\extd}{\extd t}$,
and $\hat x$ is identified with $\frac{\extd}{\extd x}$. Using this,
it can be
explicitly checked that (\ref{kjhsdvc}) gives the deformed
multiplication for this algebra up to $O(\hbar^3)$. The third order
part of (\ref{ujzvsc}) can be given by
\begin{eqnarray*}
G^{(3)}\,=\,(e_i e^j e^k\tens e_ke_j e^i-
e_k e^j e_i\tens  e^i e_j e^k-2\, e^i e_j
e^k\tens
e_i e^j e_k)/96\ ,
\end{eqnarray*}
where we sum over $i,j,k$. Note, however, that this expression
is not unique in the same manner that (\ref{kjhsdvc}) at order $\hbar^2$ is not unique as we have seen. With more work
one may exploit the non-uniqueness and expect to achieve the features in Section 5.4 with respect to  the  coproduct $\Delta_F$ as well. Note that we do not necessarily expect a unique $F,F^{-1}$ for any given Lie algebra (the uniqueness proposed in Section~5.4 was for a universal $F,F^{-1}$ applicable to all Lie algebras). 

\section{Mackey quantisations $C^\infty(N)\lcross U_\hbar(g)$ of
Homogeneous spaces as cochain twists}

Suppose that a Lie group $G$ with Lie algabra ${\frak g}$
  acts on a manifold $N$. In this case there is a
standard 'quantisation' for the system due to Mackey and much used in
physics, in which the initial algebra is $C^\infty(N)\tens S(\g)\subset C^\infty(N\times \g^*)$ (i.e. functions polynomial in the $\g^*$ direction). This is deformed or quantised to the cross product
$C^\infty(N)\lcross U_\hbar(\g)$.
Here, for $v\in  {\frak g}$ and $f\in C^\infty(N)$,
$(v\la_\hbar f)(x)=\hbar\, f'(x;v(x))$, and $U_\hbar({\frak g})$ has
the relation
$vw-wv=\hbar\,[v,w]$ in terms of the Lie bracket $[\ ,\ ]$ on ${\frak
g}$.  This algebra acts on the $L^2$ sections of a bundle whose fiber over
$x\in N$
is a representation of the stabiliser of $x$ in $G$.
In this section we show that the results of Section \ref{jkascb} may be
extended to this case also. The theory here reduces to that of
Section \ref{jkascb} when $N$ is a point. 

Note that $M=N\times \g^*$ is indeed a Poisson manifold, because the quantisation above
can be viewed as a flat deformation. Its Poisson bracket has a semidirect product form
\[ \{f,g\}=0,\quad \{v,f\}=v\la f,\quad \{v,w\}=[v,w]\]
for $f,g\in C^\infty(N)$ and $v,w\in\g^*$. 
Our goal is to lift this Poisson bivector to an element of a suitable $\CL\tens\CL$ and hence to 
a cochain $F$ at least to order $O(\hbar^2)$, i.e. to express the Mackey quantisation
as a cochain twist.

\subsection{To first order}  We will be extending the results from the previous `CBH' case in Section~5; we denote the  cochain components there  by  $G^{(i)}_{CBH}$.

\begin{defin}
Take a dual basis $(e_i,e^i)$ with $e_i\in {\frak g}$ and $e^i\in
{\frak g}^*$. Then define
some vector fields on $M=N\times \g^*$ by the following formulae, where $v\in
{\frak g}\subset S(\cg)$ and $g\in C^\infty(N)$.
\begin{eqnarray*}
\check e_i( v) \,=\, {[e_i,v]}\  &,&
\check e_i( g) \,=\,   0\ ,  \cr
\check e^i( v) \,=\, e^i(v)\  &,&
\check e^i( g) \,=\,   0\ ,  \cr
\check c_i( v) \,=\, 0\  &,&
\check c_i( g) \,=\,   e_i\la g\ .
\end{eqnarray*}
\end{defin}

\begin{propos}\label{ikcyuv}
The Poisson structure described is given by
\begin{eqnarray*}
G^{(1)} \,=\, -(\check c_i+\check e_i/2)\tens \check e^i\ .
\end{eqnarray*}
\end{propos}
\noindent {\bf Proof:}\quad Recall that for $G^{(1)} =\sum X\tens Y$
 we have $\{a,b\} \,=\, \sum X(a)\, Y(b)-Y(a)\,X(b)$. Hence
 for $v,w\in {\frak g}$ and $g,k\in C^\infty(N)$:
\begin{eqnarray*}
\check e_i( v)\, \check e^i( w) -\check e^i( v)\, \check
e_i( w)
  &=& {[e_i,v]}\, e^i(w) -
\, e^i(v) \, {[e_i,w]} \,=\, {[w,v]} -
{[v,w]}\,=\,
2\,{[w,v]}\ , \cr
\check e_i( v)\, \check e^i( g) - \check e^i( v)\,
\check e_i( g) &=& 0\ , \cr
\check e_i( g)\, \check e^i( k)  -
\check e^i( g)\, \check e_i( k)&=& 0\ , \cr
\check c_i( v)\, \check e^i( w) -\check e^i( v)\, \check
c_i( w)
  &=& 0\ , \cr
\check c_i( v)\, \check e^i( g) - \check e^i( v)\,
\check c_i( g) &=&
- e^i( v)\,  (e_i\la g)\,=\, -  v\la g\ , \cr
\check c_i( g)\, \check e^i( k)  -
\check e^i( g)\, \check c_i( k)&=& 0
\end{eqnarray*}
as desired. \eproof 

\begin{propos}
The Lie brackets between the $\check e_i$, $\check c_i$ and $\check
e^i$ are
given as follows. The $\check c_i$ commute with both the
$\check e_j$ and the $\check e^j$. The $\check e_i$ have the usual
Lie bracket
for ${\frak g}$ among themselves. The $\check e^i$ commute among
themselves.
The $\check c_i$ have the usual Lie bracket
for ${\frak g}$ among themselves. The bracket of the $e_i$ with the
$e^j$ is given by the coadjoint action. Thus the given fields form a Lie algebra 
 $\CL=\g\rcross \g^*\oplus \g$.
\end{propos}
\noindent {\bf Proof:}\quad Check against $ v$ and $g$. The
most difficult ones are:
\begin{eqnarray*}
[\check c_i,\check c_j]( g) &=& (e_i\la(e_j\la
g)-e_i\la(e_j\la g))
\,=\, ([e_i,e_j]\la g)\ ,\cr
[\check e_i,\check e^j]( v) &=& - \check
e^j({[e_i,v]})\,=\,
-e^j([e_i,v])\ .
\end{eqnarray*}
The last equation shows the coadjoint action, ${\rm
coad}_w(\psi)=-\psi\circ\ad_w$.
\eproof

We shall use this Lie algebra to induce the Mackey quantisation. We have already seen above that this is sufficient at order $\hbar$.

\begin{propos} The preconnection is for $v,w\in\cg$ and
$f,g\in C^\infty(N)$:
\begin{eqnarray*}
\varnabla_{ v} \extd  w &=& \extd{[v,w]}/2\ ,\cr
\varnabla_{v} \extd g &=& \extd(v\la g)\ ,\cr
\varnabla_{f} \extd  w &=& 0\ ,\cr
\varnabla_{f} \extd  g &=& 0\ .
\end{eqnarray*}
\end{propos}
\noindent {\bf Proof:}\quad We use the formula involving Lie derivatives;
\begin{eqnarray*}
\varnabla_{a}\xi &=& \check e^i(a)\,{\mathcal L}_{\check c_i+\check
e_i/2}\xi
-(\check c_i+\check e_i/2)(a)\, {\mathcal L}_{\check e^i}\xi\ .
\end{eqnarray*}
This gives
\begin{eqnarray*}
\varnabla_{ v}\extd  w &=&
e^i(v)\,\extd(\check c_i+\check e_i/2)( w)
-{[e_i,v]}\,\extd(\check e^i)( w)/2 \cr
&=&   e^i(v)\,\extd {[e_i,w]}/2 \ ,  \cr
\varnabla_{ v}\extd  g &=&
e^i(v)\,\extd(\check c_i+\check e_i/2)( g)
-{[e_i,v]}\,\extd(\check e^i)( g)/2 \cr
&=& e^i(v)\,\extd (e_i\la g)\ ,\cr
\varnabla_{f} \extd  w &=&
  -(e_i\la f)\,\extd(\check e_i)( w)\,=\,0 \ ,\cr
\varnabla_{f} \extd  g &=&
  -(e_i\la f)\,\extd(\check e_i)( g)\,=\,0\
\end{eqnarray*}
as required. \eproof

\begin{propos}
The curvature is given by, for $v,w,z\in{\frak g}$ and $f,g,h\in
C^\infty(N)$:
\begin{eqnarray*}
R( v, w)(\extd z) &=&
-\extd({[[v,w],z]})/4 \ ,\cr
R( v, w)(\extd g) &=& 0\ , \cr
R(v, g)(\extd z) &=&
R(v, g)(\extd h) \,=\, 0\ , \cr
R( f, g)(\extd z) &=&
R(f,g)(\extd h) \,=\, 0\ . \cr
\end{eqnarray*}
\end{propos}
\noindent {\bf Proof:}\quad First we need to state the Lie brackets
of the
vector fields, using $[\hat a,\hat b]=\widehat{\{a,b\}}$:
\begin{eqnarray*}
[\hat{ v},\hat{ w}] &=& \widehat{\{ v, w\}}\,=\,
\widehat{{[v,w]}}\ ,\cr
[\hat{ v},\widehat{ g}] &=& \widehat{\{ v, g\}}
\,=\, \widehat{(v\la g)}\ ,\cr
[\widehat{ f},\widehat{ g}] &=& \widehat{\{ f,
g\}}\,=\,0\ .
\end{eqnarray*}
Now we find the curvatures.   For $R(v,w)(\extd z)$ the computation is as in Proposition~5.1.1 with the same result. For the new case we have 
\[ R( v, w)(\extd g) = 
\extd(v\la(w\la g)-w\la(v\la g)-[v,w]\la g)=0\]
since $\la$ is a representation of the Lie algebra (if it were a cocycle representation we would
have curvature). In the expressions for $R(v, g)$ and
$R( f, g) $
every term is individually zero.\eproof

\subsection{To second order} The multiplication on
$C^\infty(N)\lcross U_\hbar(
\mathfrak{g})$
is given by
\begin{eqnarray*}
(f\tens\underline v)(g\tens \underline w) &=& f.(\underline
v_{(1)}\la_\hbar g)
\tens \underline v_{(2)}\bullet\underline w\ ,
\end{eqnarray*}
where the second product is in $U_\hbar(g)$ and the coproduct
is the usual $\Delta v=v\tens 1+1\tens v$ for $v\in \mathfrak{g}$. In
terms of
deformations of $C^\infty(N)\tens S(\mathfrak{g})$ we can decompose
the
order $\hbar^2$ part of the product as
\begin{eqnarray}\label{jkhvj}
&& f.g\tens(\hbar^2\ {\rm part\ of}\ \underline v\bullet\underline w)
+
f.(\hbar\ {\rm part\ of}\ \underline v_{(1)}\la_\hbar g) \tens
(\hbar\ {\rm part\ of}\ \underline v_{(2)}\bullet\underline w) \cr
&& +\, f.(\hbar^2\ {\rm part\ of}\ \underline v_{(1)}\la_\hbar g)
\tens
(\hbar^0\ {\rm part\ of}\ \underline v_{(2)}\bullet\underline w)\ .
\end{eqnarray}
The first term of (\ref{jkhvj}) is given by $\hbar^2\, (1\tens G^{(2)}_{CBH}{}_1)\tens (1\tens G^{(2)}_{CBH}{}_2)$, where  the final suffices $1,2$ denote the two pieces of  $G_{CBH}^{(2)}$ (summation understood). 

The second term of (\ref{jkhvj}) is, where hat denotes ommission,
\begin{eqnarray*}
&& \hbar \,  f(x).\sum_{i} (v_i\la g)(x)\tens (\hbar\ {\rm part\
of}\  (v_1\dots\hat v_i
\dots v_n)\bullet (w_1\dots w_m  )) \ .
\end{eqnarray*}
We can separate this into two stages, first the moving the $v_i$
stage, and then the
$\bullet$ multiplication. The first is given by $\hbar\,(1\tens
\check e^i)\tens
(\check e_i\tens 1)$, and the second by $\hbar\,  (1\tens G^{(2)}_{CBH}{}_1)\tens (1\tens G^{(2)}_{CBH}{}_2)$ as above.
The third term of (\ref{jkhvj}) is,
\begin{eqnarray*}
\hbar^2\, f.\sum_{i<j} (v_iv_j\la g)\tens v_1\dots\hat v_i\dots \hat
v_j
\dots v_nw_1\dots w_m \ .
\end{eqnarray*}
This is given by $\frac12 \hbar^2(1\tens \check e^i \check e^j)\tens
(1\tens \check e_i \check e_j)$, giving in total
\begin{eqnarray*}
G^{(2)}\,=\,  (1\tens G^{(2)}_{CBH}{}_1)\tens (1\tens G^{(2)}_{CBH}{}_2) + \sum_i (1\tens G^{(1)}_{CBH}{}_1 \check e^i)\tens
(\check e_i\tens G^{(1)}_{CBH}{}_2) + \frac12\sum_{ij}
  (1\tens \check e^i \check e^j)\tens
(1\tens \check e_i \check e_j)
\end{eqnarray*}

A special case is of course $N=G$ and action by left translation. Then the Mackey quantisation $C^\infty(G)\lcross U_\hbar(\cg)$ is a quantisation of $T^*G=G\times\cg^*$ and the Poisson-bracket above becomes the standard sympletic structure on $T^*G$. In that case we have an actual connection $\nabla$ in Section~6.1. 

 \subsection{Special case of $T^*G$}
In general we have a Poisson map $T^*N\to N\times \g^*$ defined
using the moment map by $(n,p)\mapsto (n,\<x_{-}(n),p\>)$ where $x_\xi$ is the vector field for the action of $\xi\in \g$ on $N$. This means a map
\[ C^\infty(N\times\cg^*)\to C^\infty(T^*N)\]
which will be surjective in the case that the action is locally transitive. In this way the Mackey quantisation results above can in principle induce quantisations of $T^*N$. 

In terms of functions on $T^*N$, $f\in C^\infty(N)$
corresponds to $\bar f=\pi^* f
\in C^\infty(T^*N)$ and $v\in {\frak g}$ corresponds to the function
$\bar v(x,p)=\<p,v(x)\>$. The elements of the algebra we order
putting all elements
of ${\frak g}$ to the right. We get the relation $\bar v\,\bar
f=\hbar\,
\overline{v\la f}+\bar f\,\bar v$. In terms of commutators,
$[\bar v,\bar f]=\hbar\, \overline{v\la f}$, and we would like this
to be given by
a Poisson bracket on $T^* N$. This means $\omega(\extd \bar
v,\extd\bar f)=
\overline{v\la f}$, or in $(x,p)$ coordinates
\begin{eqnarray*}
\frac{\partial \bar f}{\partial x^j}\,v(x)^j &=&
\omega^{(i+n)j}\,\frac{\partial \bar v}{\partial p_i}\,\frac{\partial
\bar f}{\partial x^j}
+\omega^{ij}\,\frac{\partial \bar v}{\partial x^i}\,\frac{\partial
\bar f}{\partial x^j} \cr
&=& \omega^{(i+n)j}\,v(x)^i\,\frac{\partial \bar f}{\partial x^j}
+\omega^{ij}\,p_k\,
\frac{\partial v^k}{\partial x^i}\,\frac{\partial \bar f}{\partial
x^j} \ .
\end{eqnarray*}
Provided the vector fields $v(x)$ for the action of $\cg$ span the tangent space at each point, this implies that $\omega^{(j+n)j}=1=-\omega^{j(j+n)}$ and all others
are zero,
i.e.\ the standard symplectic form on $T^*N$. Now we calculate
\begin{eqnarray*}
\omega(\extd \bar v,\extd\bar w) &=& v^j\,p_k\,\frac{\partial
w^k}{\partial x^j}
-w^j\,p_k\,\frac{\partial v^k}{\partial x^j} \,=\, \overline{[v,w]}\ .
\end{eqnarray*}

These conditions are met for  $T^*G$ with Mackey quantisation $C^\infty(G)\lcross U_\hbar(\g)$ where the action is by left translation.  Hence in this case the formulae in Section~6.1 define an actual compatible connection $\nabla$ on $T^*G$. This model will be investigated further elsewhere. 

\section{Quantum groups $\C_q[G]$ and related examples}

Here, for completeness, we show that Drinfelds original construction of quantum groups
$\C_q[G]$ (this means more precisely the dual of Drinfeld's construction) can also be formulated
as a cochain module twist. This is not fundamentally new but a useful point of view that motivated the above.  Actually, this is a general observation for any twist, in the setting of Section~2. Thus, if we are 
are interested in $A$ an initial undeformed Hopf algebra, let $A'$ be a dually paired Hopf algebra  and  $H=A'\tens A'{}^{\rm op}$. The use of a dual here is of course avoided if we work with comodule twists rather than module ones. In our case $A'$ acts on $A$ from the left by $h\la a=(\id\tens h)(\Delta a)$ and $A'{}^{\rm op}$ acts by $h\la a=(h\tens\id)(\Delta a)$. In this way $A$ becomes $H$-covariant as an algebra. Now let $f\in A'\tens A'$ be a cochain. This induces a cochain
\[ F=f_{13}f^{-1}_{24}\in H\tens H \]
where the suffices refer to the position in the four-fold tensor power of $A'$. The modified algebra $A_F$ induced by this cochain has product
\[ a\bullet b=\cdot(f_{24}f^{-1}_{13}\la (a\tens b))=f(a_{(1)},b_{(1)})a_{(2)}b_{(2)}f^{-1}(a_{(3)},b_{(3)})\]
which from the Drinfeld point of view is a usual (co)twist of the Hopf algebra $A$ into a coquasiHopf algebra.  Moreover, $A_F$ will be covariant under $H_F$ with coproduct
\[ \Delta_F(h\tens g)=f_{13}f^{-1}_{24}\Delta_{13}(h)\Delta_{24}(g)f_{24}f^{-1}_{13}\]
where the products are in the square of  $A'\tens A'{}^{\rm op}$. If we denote by $A'_f$ the usual twist of $A'$ by $f$ with coproduct $\Delta_f(h)=f(\Delta h)f^{-1}$. We see that $H_F=A'_f\tens (A'_f)^{\rm op}$, all as potentially quasi-Hopf algebras. 
This is the general situation and amounts to a recasting of the standard Drinfeld (co)twist of $A$ as a module-algebra cochain twist for a suitably doubled up $H$. 

We can of course apply this to Drinfelds example where $A=\C[G]$ is a suitable form of the coordinate ring of a classical simple Lie group with Lie algebra $\cg$, and $A'=U(\cg)$. It is understood that we  extend the above to allow formal powerseries in a parameter $\hbar$ (with $q=e^{\hbar\over 2}$).  In this case Drinfeld showed the existence (as a formal power-series) of a suitable $f$ such that $A'_f\isom U_q(\cg)$ as a Hopf algebra. Here the coboundary $\phi_f$ and hence the associator 
$\Phi_f$ are nontrivial but $\phi_f$ is central in the sense that $A'_f$ remains an ordinary Hopf algebra. This translates in the above reformulation into the statement that $F$ and its coboundary $\phi_F\in H\tens H\tens H$ are nontrivial but that $\phi_F$ acts trivially on $A$, so that $A_F\isom \C_q[G]$ remains associative, namely the usual quantum group coordinate algebra. We have $\CL=\cg\oplus\cg^{\rm op}$ as the incuding Lie algebra from our current point of view. We can further construct $\Omega(\C_q[G])=\Omega(G)_F$ and compute its associated preconnection and curvature (which is nonzero). These results are the same as those presented in \cite{BegMa:sem} (albeit from a different supercoquasiHopf algebra point of view) so we do not repeat them here.

Note that Drinfeld's $f$ is not known very explicitly (except at lowest order where the Poisson-bracket induced by the above is the usual Drinfeld-Sklyanin one), however its existence holds very generally and in a canonical way for quantum-group related examples. Moreover, it can happen that $H=U_q(\cg)$ and $F=f$ may again have $\phi_F$ acting trivially on a particular classical algebra, such as on a highest weight orbit.  Here $\CL=\cg$ is the inducing Lie algebra and we use Drinfeld's cochain without any doubling. This was the case in \cite{DGM:mat} where it was an associative quantum sphere was constructed in this way with $\CL=su_2$. This is therefore an early example of the cochain-quantisation method genuinely used.  

\section{Hidden nonassociativity}

In the above we have recovered, at least to some order, several
standard associative quantum algebras of interest in physics as cochain
twists (we do not just  mean q-deformed or quantum group examples). Here
the cochains are not required to be cocycles and this relaxation appears to be necessary. It means, however, that even though the algebra of 'functions' happens to remain
 associative, there is
an underlying nonassociativity behind the scenes in all these quantum algebras.
We now turn to this aspect.

First of all, as our algebras become quantized, their covariance Lie algebra $\CL$ gets deformed to
a quasi-quantum group $U(\CL)_F$ as explained in Section~2. These are in principle 'noncoassociative' and are looked at for our various examples in 
Section 8.1. Next, our  quantum algebras are all equivalent in a certain monoidal categorical sense
to the unquantised algebras, with the result that not only the
algebras but all functorial constructions on them are similarly
quantised, for example differential forms on the classical phase space and the Dirac operator deform naturally to the quantum algebras, but nonassociatively. We
consider these in Sections 8.2 and 8.3 respectively.

\subsection{The quasiHopf algebras $U(\CL)_F$}

From Section~2, the deformed algebra $A_F$ remains covariant, but under the quasi-Hopf algebra $H_F$. In our cases of interest $H=U(\CL)$ which is also the algebra of $H_F$. Its coproduct, however, is modified to
\[ \Delta_F(X)=F(X\tens 1+1\tens X)F^{-1}=X\tens 1+1\tens X-\hbar[G^{(1)},X\tens 1+1\tens X]+O(\hbar^2)\]
for any $X\in \CL$. The leading order data here defines a quasi-Lie bialgebra $(\CL,\delta,\psi)$ 
where 
\[ \delta X=[X\tens 1+1\tens X,G^{(1)}]=\ad_X(G^{(1)})\in \CL\tens\CL.\]
If it happens that $\delta$ obeys the cojacobi identity
\[ (\delta\tens\id)\delta X+{\rm cyclic}=0\]
then we have an ordinary Lie bialgebra, which means that at lowest order at least, $U(\CL)_F$ remains an ordinary (not quasi) Hopf algebra. This is already the case for the Drinfeld twist examples in Section~7 and indeed the $O(\hbar^2)$ part $\psi$ of $\phi$ is a multiple of the ad-invariant Cartan tensor $n\in\Lambda^3(\cg)$ defined by the Killing form.  In general the cojacobiator above is given by $\ad_X(\psi)$ and this is the fundamental reason why the covariance algebra remains (co)associative for such examples. But let us see how the situation fares for our non-quantum group examples.

 Thus, $\psi$ is given for the $\CL=so(1,3)$ example in Section~4.5 and from the expression there one may readily compute that
 \[ \ad_{X_i}(\psi)\ne 0,\quad \ad_{Y_i}(\psi)=0\]
 i.e. rotationally invariant (as to be expected as the whole construction is) but not invariant under boosts.  Thus our `sphere at infinity' example in Section~4 gives us a quasi-Hopf algebra version of $\CL=so(1,3)$.  
 
 Next up, we CBH or $U_\hbar(g)$ example in Section~5, we have seen in Section~5.4 what the twisted coproduct $\Delta_F$ looks like to $O(\hbar^3)$ when acting on $\C[\cg]\subset H=U(\cg\rcross\cg^*)=U(\cg)\rcross\C[\cg]$. We have seen that this twists to a local form of the classical coordinate ring $\C[G]$. On the other hand the coproduct of $U(\cg)$ remains unchanged after twisting to this order because all elements are $g$-invariant under commutator in the bigger algebra and hence
 \[ [\Delta v,G^{(i)}]= [v\tens 1+1\tens v,G^{(i)}]=[v,G^{(i)}_1]\tens G^{(i)}_2+G^{(i)}_1\tens [v,G^{(i)}_2]=0.\]
 This is clear since the expressions involve only paired bases and dual basis, the Killing form etc. and the commutators of $v\in\cg$ are given by the adjoint and coadjoint actions.  Such invariance would be a reasonable requirement to all orders for any universal formula for  $F,F^{-1}$. 
 Since $H=U(\cg)\rcross\C[\cg]$ is generated by
 $U(\cg)$, $\C[\cg]$ and has the same algebra after twisting, we conclude that 
\[ U(g\rcross g^*)_F=U(g)\rcross\C[g]_F\]
as an ordinary Hopf algebra. Moreover, this is locally isomorphic to $U(g)\rcross \C[G]=D(U(g))$, the Drinfeld quantum double, which is an ordinary Hopf algebra. It is known that $U_\hbar(\cg)$ is always covariant under $D(U(\cg))$ and this was explored for $U_\hbar(su_2)$ (the so-called universal `fuzzy sphere') in \cite{BatMa:non}.  In this case the background covariance becomes a quantum group covariance but remains associative (there is still hidden nonassociativity, see below).

Finally, the Mackey case is a nontrivial extension of the CBH case and has a larger symmetry group. We do not make the full analysis here but suffice it to say that one reason that the CBH case works from an algebraic point of view to the fact that $U_\hbar(\cg)$ is a (cocommutative) Hopf algebra and therefore has a larger quantum group covariance based on the Drinfeld double; there is no such argument for the general Mackey case but there are cases when $N$ is itself a group and acts back on $G$ such that the Mackey quantisation in an algebraic form becomes a bicrossproduct Hopf algebra, see \cite{Ma:book}. In such cases one might expect similar behaviour to the CBH case above, but not in general. The bicrossproduct case includes the deformed Poincar\'e group for the noncommutative spacetimes mentioned in Section~5.6, see \cite{MaRue:bic}.

 \subsection{Quasiassociative quantum differential calculi}

Next, by applying the same cochain twist to the classical exterior algebra, we obtain
  noncommutative differential calculi on our various quantisations. We mean differential calculus  in the sense of noncommutative
geometry but in a monoidal weakly associative category, and we will
see that our calculi on these examples are indeed nonassociative.
They do, however, have the merit of classical dimensions in each
degree. Again, the quantum group case in Section 7 was already covered in an equivalent form 
in \cite{BegMa:sem} with the main result  that the resulting $\Omega(\C_q[G])$ have curvature and hence are not associative even though the quantised algebra happens to be. But let us see how our non-quantum group-related examples fare.  

First, for our sphere at infinity example. The given vector fields act by Lie derivative,
which commutes with the $\extd$ operator. Some large calculations give
the special cases:
\begin{eqnarray*}
x\bullet\extd x &=& x.\extd x+\frac{\hbar}{2\,z}
\big(x^2 y.\extd x+x(1-x^2).\extd y\big)
+\frac{\hbar^2}{8}\,x.\extd x +O(\hbar^3)\ , \cr
x\bullet\extd y &=& x.\extd y+\frac{\hbar}{2\,z}
\big(x y^2.\extd x+y(1-x^2).\extd y\big)
+\frac{\hbar^2}{8}\,x.\extd y +O(\hbar^3)\ , \cr
y\bullet\extd x &=& y.\extd x-\frac{\hbar}{2\,z}
\big(x (1-y^2).\extd x+x^2 y.\extd y\big)
+\frac{\hbar^2}{8}\,y.\extd x +O(\hbar^3)\ , \cr
y\bullet\extd y &=& y.\extd y-\frac{\hbar}{2\,z}
\big(y (1-y^2).\extd x+x y^2.\extd y\big)
+\frac{\hbar^2}{8}\,y.\extd y +O(\hbar^3)\ .
\end{eqnarray*}
Moreover, since the connection $\nabla$ arising from the noncommutativity of the calculus at lowest degree turned out to be the Levi-Civita one, and since this has (constant) curvature, we know that the exterior algebra of this quantised sphere is necessarily nonassociative.

Next, for the CBH or  $U_\hbar(g)$ example, again we have found expressions for the curvature in terms of a double commutators. Whether or  not this vanishes depends on the Lie algebra in question: for the Heisenberg Lie algebra for example, one has $R=0$. However, for a simple Lie algebra  such double commutators will not vanish and there is curvature, hence nonassociativity of $\Omega(U_\hbar(\cg))$.  One can also
write the deformed calculus explicitly:
\[ v\bullet\extd w= v\extd w+\hbar\beta\, (e^i\la v)\extd e_i\la w)+ O(\hbar^3)=v\extd w+\hbar\beta\, \extd[v,w]+O(\hbar^3)\]
\[ \extd w\bullet v=(\extd w)v+\hbar\alpha\, (\extd e_i\la w)e^i\la v+O(\hbar^3)=(\extd w)v+\hbar\alpha\, \extd[v,w]+O(\hbar^3)\]
where we use the expression for $F^{-1}$ in Section~5.4. The second order terms fail to contribute because $e^i\la v=v^i.1$ is degree zero (it acts by differentiation) which is then killed by $\extd$. Therefore the only term that could contribute in the first line, for example, is from $Q_2^2$ i.e. $(e^ie^j\la v)\extd e_ie_j\la w$ which is zero because of the second differentiation on $v$. The difference between the two expressions is of course 
\[ v\bullet\extd w-\extd w\bullet v={\hbar \over 2}\extd [v,w]+O(\hbar^3),\]
 i.e. the Poisson-compatible preconnection as it should be. The example of  `noncommutative spacetime' with nonassociative differentials
\[ t\bullet \extd x_i-\extd x_i\bullet t=\extd t\bullet x_i-x_i\bullet\extd t={\hbar\over 2} \extd x_i\]
is very different from the associative (but not canonical) differential calculus usually used for this model. It represents a different approach that may overcome some of the structural problems encountered previously (such as to find the canonical Dirac operator, see below). Note that in the physical application the deformation parameter that we have denoted $\hbar$ should be denoted by another symbol and is expected to be of the order of the Planck time $\sim 10^{-44}$s.

Finally, the Mackey quantisaition case is more complicated but from the curvature computations in Section~5 we conclude again that the natural deformed calculus $\Omega(C^\infty(N)\lcross U_\hbar(\cg))$ is nonassociative at least for simple Lie algebras, because the CBH part is. It is interesting to note that the curvature comes from this part alone.

\subsection{Isospectral quantum Dirac operator}

Here we conclude with an example to demonstrate that the categorical deformation
method outlined in Section~2 is very powerful indeed and quantizes almost any natural construction. 
In other words, when a quantisation is expressed as a cochain module algebra twist this has great consequences.

Specifically, in another approach to noncommutative geometry it is normal to look for an analogue of the
Dirac operator in the form of a 'spectral triple' \cite{Con} obeying some axioms. These axioms 
are natural from an associative point of view but it is well known that important examples
such as $\C_q[G]$ do not admit operators obeying exactly those axioms. We see by contrast
that there is a natural deformation of any classical Dirac operator on the classical phase space but it will
obey a variation of Connes axioms due to the hidden nonassociavity in the underlying differential calculus and elsewhere. We explain now that such an approach
agrees with recent 'isospectral deformation' proposal for the Dirac operator on $\C_q[SU_2]$ 
in \cite{Landi}. On the other hand, it is more categorical and works in principle for all quantum groups $\C_q[G]$, and moreover works for
our more conventional quantisations such as $U_\hbar(g)$ and the Mackey quantisation to provide (in principle at least) some type of Dirac operators on them.

  We consider for the sake of discussion only the case where the classical and hence quantum cotangent and spin bundles are trivial so that  the spin bundle in particular has the form $V\tens A$ where $A=C^\infty(M)$ and $V$ is ostensibly the representation space for the spin group. We do require everything to be covariant under a background Lie algebra $\CL$ (or Hopf algebra $H$) to induce the quantisation given a cochain. This is not a problem in the case $M=G$ a Lie group  (a covariant Dirac operator). The short version of the quantisation is then as follows: we consider the classical Dirac operator $D:V\tens A\to V\tens A$ and to this we apply the functor $\CT$ in Section~2 to obtain a map $\CT(D):\CT(V\tens A)\to \CT(V\tens A)$. 
As in Section~2 we have to allow that although $\CT$ acts as the identity on objects and morphisms (so $\CT(A)=A$ which becomes the deformed algebra $A_F$, $\CT(V)=V$, $\CT(D)=D$), it is nontrivial as a monoidal functor and in the sense of potentially nontrivial natural isomorphisms $\CT(V)\tens\CT(A)\isom \CT(V\tens A)$ with certain properties in relation to $\tens$. We refer to \cite{Ma:book} for an introduction. Here these isomorphisms are given by the action of $F^{-1}$. Putting these facts together, we have the deformed Dirac operator:
\[ D_\bullet:A_F\tens V\to A_F\tens V,\quad D_\bullet(a\tens v)=F\la D(F^{-1}\la (a\tens v)).\]
The main thing to note about this construction is that since it is given by conjugation by $F$ as an operator, it does not change the spectrum in the  Hilbert space setting. It should be remarked that one would still need a lot of analysis to make these remarks fully precise. 

The only other subtlety is to identify $A_F$ (which is the same vector space as $A$ with the deformed product) explicitly as $\C_q[G]$ in the case $A=\C[G]$. This is not trivial but we note that 
 $\CT$ respects sums so if one has made a Peter-Weyl decomposition of $A$ into a  direct sum of matrix algebras (as is possible for $A=\C[G]$ for simple $G$) and likewise decompose $\C_q[G]$ in its Peter-Weyl decomposition, we can identify the summands as matrix coalgebras.  This is our interpretation of the proposed Dirac operator in \cite{Landi}.  On the other hand, $D_\bullet$ lives in a nontrivial monoidal category and has properties in which the non-associativity of the category will surely play a role. It is known that the axioms in \cite{Con} are not satisfied and we would propose to replace them by ones that take this hidden nonassociavity into account.
 
 Finally, the longer answer to the deformed Dirac operator here is to 'get inside' its construction. One can do this too in principle as we outline now. Thus, we break $D$ into a series of morphisms all covariant under our background Hopf algebra.  We also suppose for the sake of discussion that $M$ is parallelizable (e.g. $M$ a Lie group) so that its (covariant) differential calculus has the form $\Omega^1(M)=A\tens\Lambda^1$ as a (trivial) bundle associated to $\Lambda^1$. We write $\extd a=(\del^i a)\tau_i$ where $\tau_i$ are a basis of $\Lambda^1$ and take this as a definition of the partial derivatives. Finally, we assume `$\gamma$-matrices' of some form $\gamma:\Lambda^1\tens V\to V$ so that 
 \[ D(a\tens v)=\del^i(a)\tens\gamma_i(v)=(\id\tens\gamma)(\extd\tens\id)(a\tens v)\]
 expresses $D$ as a composition of morphisms (here $\gamma_i=\gamma(\tau_i\tens (\ ))$ would be the more conventional point of view). Note that the defining relations among the $\gamma$ also needs to be invariant under the background symmetry, e.g. by an invariant metric.

 We now define $\gamma_\bullet(\tau \tens v)=\gamma(F^{-1}\la (\tau\tens v))$ for all $\tau\in \Lambda^1$, $v\in V$ by the same reasoning as above, i.e. the functor $\CT$. Likewise we have $\extd_\bullet=F\la\extd$ when $\extd:A\to A\tens\Lambda^1$. Note that in the above we have not deformed $\extd$ and indeed this is not deformed if we consider it to $\Omega^1$ and identify $\CT(\Omega^1)=\Omega^1$ in the deformed theory, $\extd_\bullet$ is a slightly different object. This now expresses the above as
 \[ D_\bullet=F D F^{-1}=(\id\tens\gamma_\bullet)\Phi_{A,\Lambda^1,V}(\extd_\bullet\tens
\id)\]
in view of the diagram:
\[\begin{matrix} \CT(A\tens V)& {  \longrightarrow}& \CT(A)\tens \CT(V)& &\\
\downarrow \extd &  &\downarrow \extd  & \quad \searrow \extd_\bullet \kern -20pt &  \\
\CT((A\tens\Lambda^1)\tens V)&\longrightarrow & \CT(A\tens\Lambda^1)\tens\CT(V)& \longrightarrow &((\CT(A)\tens \CT(\Lambda^1))\tens\CT(V)\\
|\, | & & & & \downarrow \Phi_{A,\Lambda^1,V}\\
\CT(A\tens(\Lambda^1\tens V))&\longrightarrow & \CT(A)\tens\CT(\Lambda^1\tens V)& \longrightarrow &(\CT(A)\tens( \CT(\Lambda^1)\tens\CT(V))\\
\downarrow \gamma &  &\downarrow \gamma &  \quad \swarrow \gamma_\bullet \kern -20pt&  \\
\CT(A\tens V)&  \longrightarrow &  \CT(A)\tens\CT(V)& 
\end{matrix}\]
The large middle cell here commutes by definition of a monoidal functor (the associator in the initial category of $H$-covariant objects is trivial). The upper and lower left cells commute because of operations on different spaces. The upper and lower right cells commute by the definitions of $\extd_\bullet$ and $\gamma_\bullet$ respectively. The horizontal arrows are all given by the action of $F$. The vertical composition on the left is $\CT(D)=D$, while the vertical composition on the right is the definition of $D_\bullet$.

As to the $\gamma_\bullet$, one should write the defining relations of $\gamma$ as commuting diagrams, apply the functor $\CT$ to obtain commuting diagrams in the deformed category, and use $F$ to interpret them as $\gamma_\bullet$ relations in a similar manner to the above. The result is, for example, a deformed set of Clifford relations involving now $\Phi$ and the deformation of the flip map induced by $F$ (as a symmetry in the category). The actual relations would depend on the classical set up which  need not be the usual classical Clifford relations if the frame group is not the usual one (e.g. one may use a Lie group to frame itself).

\end{document}